\documentclass[reqno]{amsart}
\usepackage{url}
\usepackage{amscd}
\usepackage{amssymb}
\usepackage{hyperref}
\usepackage[all]{xy}

\renewcommand{\theenumi}{(\roman{enumi})}

\newtheorem{thm}{Theorem}[section]
\newtheorem{prop}[thm]{Proposition}
\newtheorem{lem}[thm]{Lemma}
\newtheorem{cor}[thm]{Corollary}

\theoremstyle{definition}
\newtheorem{rem}[thm]{Remark}
\newtheorem{defn}[thm]{Definition}
\newtheorem{exa}[thm]{Example}

\renewcommand{\theenumi}{(\roman{enumi})}
\numberwithin{equation}{section}

\def\Z{\mathbb{Z}}
\def\Q{\mathbb{Q}}
\def\F{\mathbb{F}}
\def\R{\mathbb{R}}
\def\C{\mathbb{C}}
\def\U{\mathbb{U}}

\def\O{\mathcal{O}}
\def\E{\mathcal{E}}
\def\A{\mathcal{A}}

\def\p{\mathfrak{p}}
\def\q{\mathfrak{q}}
\def\P{\mathfrak{P}}

\def\a{\mathfrak{a}}

\def\Qp{\Q_p}
\def\Fp{\F_p}
\def\Zl{\Z_\ell}
\def\ad#1{\mathbf{A}_{#1}}
\def\id#1{\ad{#1}^\times}
\def\ab{\mathrm{ab}}
\def\Qab{\Q^{\ab}}
\def\Kab{K^{\ab}}
\def\M{\mathrm{M}}
\def\GL{\mathrm{GL}}
\def\SL{\mathrm{SL}}
\def\Fr{\mathrm{Fr}}
\def\Tr{\mathrm{Tr}}
\def\N{\mathrm{N}}
\def\Gal{\mathrm{Gal}}
\def\End{\mathrm{End}}
\def\Aut{\mathrm{Aut}}
\def\ord{\mathrm{ord}}
\def\jd{j}
\def\gammatwod{\gamma_2}
\def\gammathreed{\gamma_3}
\def\Ec{\E_{\mathrm{can}}}

\def\H{\mathfrak{H}}
\def\bmu{\boldsymbol{\mu}}
\def\leg#1#2{{\bigl(\textstyle\frac{#1}{#2}\bigr)}}
\def\legt#1#2{{\bigl(\textstyle\frac{#1}{#2}\bigr)_{\hskip-1pt 2}}}
\def\legf#1#2{{\bigl(\textstyle\frac{#1}{#2}\bigr)_{\hskip-1pt 4}}}

\def\legtF#1#2#3{{\bigl(\textstyle\frac{#1}{#2}\bigr)_{\hskip-1pt 2,#3}}}
\def\legfF#1#2#3{\bigl(\textstyle\frac{#1}{#2}\bigr)_{\hskip-1pt 4,#3}}
\def\legtwo#1#2{\bigl(\textstyle\frac{#1}{#2}\bigr)_{\hskip-1pt 2}}
\def\legfour#1#2{\bigl(\textstyle\frac{#1}{#2}\bigr)_{\hskip-1pt 4}}
\def\hilb#1#2{[#1,#2]}

\def\hilbt#1#2#3{\hilb{#1}{#2}_{2,#3}}

\def\isom{\,\mbox{$\xrightarrow{\hskip6pt}\hskip-13pt\raisebox{3.5pt}{\footnotesize$\sim$}\hskip6pt$}}
\def\tinymatrix#1#2#3#4{\bigl(\begin{smallmatrix}#1 & #2 \\ #3 & #4\end{smallmatrix}\bigr)}

\def\abcd{\tinymatrix{\alpha}{\beta}{\gamma}{\delta}}
\def\shortoverline#1{\hskip 2pt\overline{\hskip -2pt #1 \hskip-2pt} \hskip 2pt}

\def\td{\tau_D}
\def\mm{\mu}
\def\eone{{\delta}_\tau}
\def\etwo{\epsilon_\tau}
\def\etd{\epsilon_{\tau_D}}
\def\etdtilde{\delta_{\tau_D}}
\def\zzero{\tau}
\def\Et#1#2{E^{(#1)}_{#2}}
\def\Ol#1{\O_{\tau,#1}}
\def\Eb{\tilde{E}}

\title[Point counting on CM elliptic curves]{Point 
counting on reductions of CM elliptic curves}

\author{K.\ Rubin}
\author{A.\ Silverberg}
\address{Mathematics Department, 
University of California, 
Irvine, CA 92697, 
USA}
\email{krubin\char`\@uci.edu}
\email{asilverb\char`\@uci.edu}

\thanks{This material is based upon work supported by the 
National Science Foundation under grants DMS-0457481 and DMS-0757807
and the National Security Agency under grants
H98230-05-1-0044 and H98230-07-1-0039.}

\begin{document}

\begin{abstract}
We give explicit formulas for the number of points on reductions of elliptic curves with
complex multiplication by any imaginary quadratic field.
We also find models for CM $\Q$-curves in certain cases.
This generalizes earlier results of Gross, Stark, and others.
\end{abstract}

\maketitle

\section{Introduction}

In this paper we give explicit formulas for the number of points on reductions of CM elliptic curves
(see Theorems \ref{mainspec} and \ref{main} and Corollary \ref{exrem}).
We also give models for CM $\Q$-curves, in certain cases (see Theorem \ref{qcurvethm}).

If $\tilde{E}$ is an elliptic curve over a finite field $\F_q$, it is 
well-known that to count the number of points in $\tilde{E}(\F_q)$,  
it suffices to determine the Frobenius endomorphism of $\tilde{E}$ over $\F_q$, 
or more precisely the trace of Frobenius acting on an appropriate 
vector space.  The best methods known for accomplishing this with 
a general elliptic curve are 
modifications of the method of Schoof \cite{schoof, schoof2} or $p$-adic methods
\cite{Satoh}.

When $\tilde{E}$ is the reduction of an elliptic curve $E$ over a number 
field $F$ with complex multiplication (CM) by an order in an imaginary quadratic 
field $K \subseteq F$, a different approach is possible.  In this case, 
as shown by Deuring \cite{deuring}, there is a Hecke character $\psi$ of $F$ 
with values in $K^\times$ such that for every prime $\P$ of 
$F$ where $E$ has good reduction, $\psi(\P) \in  K = \End(E) \otimes \Q$ 
reduces to the Frobenius endomorphism of $E$ modulo $\P$.  Thus if one 
can compute the Hecke character $\psi$, one can determine the number 
of points on every reduction of $E$, including the original curve $\tilde{E}$. 
If $\tilde{E}$ is an ordinary elliptic curve over $\F_q$, then $\tilde{E}$ 
is always the reduction modulo $\P$ of some CM elliptic curve 
$E$ defined over some number field $F$.  
The field $K$ 
determines the Frobenius endomorphism of $\tilde{E}$ over $\F_q$ up to a 
root of unity in $K$ (generally $\pm1$).  The computation of the Hecke character 
of $E$ can be viewed as the determination of this 
root of unity, for every prime $\P$ of $F$.

This CM approach has been carried out in special cases by several authors.  
The Hecke character of $E$ was computed by Gross \cite{gross,grossbk}
when $\End(E)$ is the maximal order in $\Q(\sqrt{-p})$ with $p$ prime and
$p \equiv 3 \pmod{4}$, and by Stark \cite{Stark} 
when $\End(E)$ is the maximal order in $\Q(\sqrt{-d})$ with squarefree 
$d \equiv 3 \pmod{4}$ and $3 \nmid d$ (i.e.,  $d\equiv 7$ or $11\pmod{12}$).
Individual special cases were done earlier by a number of people, dating
back to Gauss;
see p.~349 of  \cite{Rajwade} for some of the
relevant references.
For further discussion of the history of this problem, see \S5 of \cite{Stark}.

In this paper we complete this program by computing,
for every imaginary quadratic field $K$, every
imaginary quadratic order $\O$, and every number field $F \supseteq K$, the 
Hecke character of every elliptic curve over $F$ with $\End(E) \cong \O$,
thereby computing the number of points on the reductions of these elliptic curves.
This extends the results of Stark and Gross to all $d$, 
including $d \equiv 1, 2\pmod{4}$ and $d \equiv 3\pmod{12}$, and to all orders, including
non-maximal orders.
Also, whenever $d \equiv 2$ or $3 \pmod{4}$,
we produce a model of a $\Q$-curve with CM by the maximal order in $\Q(\sqrt{-d})$.
(There are no $\Q$-curves with CM by the maximal order in $\Q(\sqrt{-d})$ 
when $d > 1$ is a product of primes congruent to $1\pmod{4}$.)

One motivation for studying this question comes from cryptography.  
For various cryptographic applications, such as finding 
``pairing-friendly" elliptic curves,
one needs to find an elliptic curve over $\Fp$ with a given number  
of points.  The usual way to do this (the ``CM method'' \cite{atkinmor}) 
produces a CM elliptic curve over a number field whose reduction $\tilde{E}/\F_p$ 
has the property that either $\tilde{E}$ or its
quadratic twist has the correct number of points.  In \cite{cmmethod} we 
use the results in this paper to give a simple efficient algorithm
for determining which of the two elliptic curves is correct.
This settles an open question of Atkin and Morain 
(Conjecture 8.1 of \cite{atkinmor}). 

We now state  our main result in the (useful) special case where $j(E) = j(\O_K)$,
with $\O_K$ the maximal order (it follows that $E$ has CM by $\O_K$). 

\begin{thm}
\label{mainspec}
Suppose $E : y^2 = x^3 + ax + b$ is an elliptic curve over a number field $F$, and 
$j(E) = j(\O_K)$ where $\O_K$ is the ring of integers of an imaginary quadratic field 
$K=\Q(\sqrt{-d}) \subseteq F$,
with squarefree $d \neq 1, 3$.
Suppose $\P \nmid 2$ is a prime of $F$
where $E$ has good reduction. 
Let $\lambda\in\O_K$ be a generator of the principal ideal $N_{F/K}(\P)$ 
and let $q = \N_{F/\Q}(\P)$.
Then
$$
\#E(\O_F/\P) = q + 1 - W\cdot \epsilon \cdot\Tr_{K/\Q}(\lambda)
$$
where
$$
W = 
\begin{cases}
\legtwo{6b\gamma_3(z_d)}{\P} & \text{if $d  \equiv 3\pmod{4}$}, \\[7pt]
\legtwo{-6bi{\gamma_3(z_d)}}{\P} & \text{if $d  \equiv 2\pmod{4}$}, \\[7pt]
\legfour{(6b)^2(j(E)-1728)}{\P} & \text{if $d  \equiv 1\pmod{4}$},
\end{cases}
$$
the $n$-th power residue symbols $\leg{c}{\P}_{n} \in \bmu_n$  
 and  the Weber function $\gamma_3$ are defined in \S\ref{nota} below,
$z_d$ is defined by 
$$
\renewcommand{\arraystretch}{1.75}
\begin{array}{|c||c|c|c|c|}
\hline
d\hskip-7pt\pmod{8}  & 2 & 3 & 6 & 7 \\
\hline
z_d & \sqrt{-d} & \frac{3+\sqrt{-d}}{2} & 3+\sqrt{-d} & \frac{-3+\sqrt{-d}}{2} \\
\hline
\end{array}
$$
and $\epsilon$ is defined by:

\medskip\noindent
$d\equiv 3\pmod{4}$: 

\medskip~\hskip-14pt
$\arraycolsep=8pt
\renewcommand{\arraystretch}{1.25}
\begin{array}{|c||c|c|}
\hline
\lambda^3~(\mathrm{mod}~4)  & 1, -\sqrt{-d} & -1, \sqrt{-d} \\
\hline
\epsilon & 1 & -1 \\
\hline
\end{array}$

\medskip\noindent
$d\equiv 2\pmod{4}$:

\medskip~\hskip-14pt
$\arraycolsep=8pt
\renewcommand{\arraystretch}{1.25}
\begin{array}{|c||c|c|}
\hline
\lambda~(\mathrm{mod}~4) & 1, -1+2\sqrt{-d}, \pm1+\sqrt{-d} 
    & -1, 1+2\sqrt{-d}, \pm1-\sqrt{-d} \\
\hline
\epsilon & 1 & -1  \\
\hline
\end{array}$

\medskip\noindent
$d\equiv 1\pmod{4}$:

\medskip~\hskip-14pt
$\arraycolsep=5pt
\renewcommand{\arraystretch}{1.25}
\begin{array}{|c||c|c|c|c|}
\hline
\lambda~(\mathrm{mod}~4)  & 1, 1+2\sqrt{-d} & 2+\sqrt{-d},\sqrt{-d} 
    & -1,-1+2\sqrt{-d} & 2-\sqrt{-d}, -\sqrt{-d} \\
\hline
\epsilon & 1 &i & -1 & -i \\
\hline
\end{array}$
\end{thm}

Our method of proof
is similar to the method of Stark \cite{Stark}, which 
follows an approach used by Rumely in his thesis and \cite{Rumely}.  
Rumely showed how to use Shimura's Reciprocity Law (for values 
of modular functions at CM points) to compute the Hecke character 
of a CM elliptic curve in certain special parametrized families.\footnote{Shimura
points out in Remark 14.12(3) of \cite{val} that 
there is a gap in Rumely's proof of Theorem 1 of \cite{Rumely}, although the 
statement of that theorem is correct in the setting of Example 1 of \cite{Rumely}.
While our method was inspired by Rumely's approach, we do not use his results.}
Rumely (Example 1 on p.~394 of \cite{Rumely})
and Stark (equation (3) on p.~1121 of \cite{Stark})
used Weber functions to write down a family 
$E_z$ of elliptic curves, parametrized 
by $z $ in the complex upper half-plane $\H$ (take $\alpha = 1$ in 
Definition \ref{esubdef} below).  When $d \equiv 3 \pmod{4}$ and 
$3 \nmid d$, then $z \in \H$ can be chosen so that $E_z$ has CM by 
the maximal order $\O_K$ of $K=\Q(\sqrt{-d})$ and 
$E_z$ is defined over the Hilbert class field $H_K$ of $K$, and in this case 
Stark computes the Hecke character of $E_z$ over $H_K$.  
If $E$ is an arbitrary elliptic curve with CM by $\O_K$ over a number field $F \supseteq K$, 
then $H_K \subseteq F$ and $E$ is isomorphic to a quadratic twist of some such 
$E_z$ over $F$, so one obtains the Hecke character of $E$ over $F$.

If either  $d$ is a multiple of $3$
or $d \not \equiv 3 \pmod{4}$, then there are 
$z \in \H$ such that $E_ z $ has CM by $\O_K$.
For all such $z $, the curve $E_ z$ is defined 
over a small but nontrivial extension of $H_K$.  
For arbitrary orders $\O$ there are 
$z \in \H$ such that $E_z$ has CM by $\O$ and $E_z$ is defined 
over a small extension $H'_\O$ of the ring class field $H_\O$ of $\O$.  
If $E$ is an elliptic curve  with CM by 
$\O$ defined over a number field $F \supseteq K$, 
then $E$ is isomorphic to some $E_z$ over $\bar{\Q}$, 
and $F$ contains $H_\O$ but $F$ need not contain $H'_\O$.  
In order to compute the Hecke character of $E$ over $F$, we need to 
determine what $\Gal(\bar{\Q}/H_\O)$ does to the torsion points of $E_z$, 
not just the action of its proper subgroup $\Gal(\bar{\Q}/H'_\O)$
on  the torsion points.  We do this in 
Proposition \ref{ptrans}, extending the Rumely-Stark method.   
This allows us to compute the Hecke characters for all elliptic curves
with  CM by $\O$ defined over $F$, for {\em every} $d$ and
$\O$ and every number field $F \supseteq K$.
Our main results are  Theorem \ref{main} and Corollary \ref{exrem},
and the heart of the proof is in Theorem \ref{ptranscor}.

In \cite{grossbk}, Gross defined a $\Q$-curve to be an elliptic 
curve that is isogenous to all of its Galois conjugates, and studied these curves 
in detail when they have CM.  In \cite{gross}, Gross exhibited equations for 
$\Q$-curves with CM by the maximal order of $\Q(\sqrt{-p})$ when $p$ is a prime 
congruent to $3 \pmod{4}$, and determined their Hecke characters.  
We use our Hecke character computations (Theorem \ref{main}) 
to exhibit equations for 
$\Q$-curves with CM by the maximal order of $\Q(\sqrt{-d})$ for all 
$d \equiv 2$ or $3 \pmod{4}$, 
and we use quadratic reciprocity over $K$ to give another expression 
(Theorem \ref{qcurvethm}) for the Hecke characters of these curves.  
When $d \equiv 3 \pmod{4}$,
the formula for the Hecke character in Theorem \ref{qcurvethm} is the one given 
by Gross (Theorem 12.2.1 of \cite{grossbk} and Proposition 3.5 of \cite{gross})
when $d$ is prime
and by Stark (Theorem 1 of \cite{Stark}) when $3 \nmid d$,
while the formula in Theorems \ref{mainspec} and \ref{main}
is of a different form.

In Example \ref{psimythex} we give a counterexample to the common myth that 
$\psi(\P)$ is necessarily in $\O$,
where $\psi$ is the Hecke character associated to an elliptic
curve with CM by an order $\O$.

The reader who wishes to avoid technical details might prefer to
start by reading the statements of Theorems \ref{mainspec},  \ref{main}, and
 \ref{qcurvethm} and Corollary \ref{exrem}, and referring back to the notation
 and supporting lemmas and propositions as necessary.

\smallskip

\noindent{\bf Outline of the paper.}
In \S\ref{nota} we introduce notation, state Shimura's 
Reciprocity Law, and describe the 
setting in which we work.  In \S\ref{CtHc} we state or work out the 
properties of the Weber functions and Dedekind's $\eta$-function that 
we need to compute Hecke characters.
In \S\ref{heckesect} (Theorem \ref{ptranscor})
we use these properties to compute 
the Hecke characters of the twists of $E_z$ mentioned above.
In \S\ref{Mt} we use Theorem \ref{ptranscor} to prove 
Theorem \ref{main} and Corollary \ref{exrem}, our main results
on Hecke characters and point counting, 
and in \S6 we compute and exhibit the
tables of values of an important function 
that appears in our formulas in Theorem \ref{main} and Corollary \ref{exrem}.
In \S\ref{Qc} we 
obtain models for 
$\Q$-curves and formulas for their Hecke characters (Theorem \ref{qcurvethm}).  
In \S\ref{yamihere} we 
give a point-counting 
result with a different flavor, under hypotheses that lead to a particularly 
simple statement.

\section{General notation}
\label{nota}

In this section we give definitions and notation that will be used in later
sections, and state Shimura's Reciprocity Theorem.

Let $\H$ denote the complex upper half-plane.
Let $i$ denote the square root of $-1$ in $\H$.
For $z\in \H$, let 
$$
L_z:= \Z+\Z z,
$$
$$
g_2(z) := 60 \sum_{0\neq\omega\in L_z} \omega^{-4} \qquad \text{and} \qquad
g_3(z) := 140 \sum_{0\neq\omega\in L_z} \omega^{-6},
$$
and  let $\wp(u;z)$ denote the Weierstrass 
$\wp$-function of $u \in \C$ for the lattice $L_{z}$.

Note that $g_k(z)$ is a modular form of weight $2k$ and level $1$, with Fourier coefficients 
in $(2 \pi i)^{2k}\Q$ (see for example \S2.2 of \cite{Shimurabook}).
Let $\eta$ denote the Dedekind eta function
$\eta(z):=e^{2\pi i z/24}\prod_{n=1}^\infty (1-e^{2\pi inz})$,
and define the Weber functions 
$$
\gamma_2(z) := 12\frac{g_2(z)}{(2 \pi i)^{4}\eta(z)^{8}} \qquad \text{and} \qquad 
\gamma_3(z) := -6^3\frac{g_3(z)}{(2 \pi i)^{6}\eta(z)^{12}}.
$$ 
Then $\eta^8$ (resp., $\eta^{12}$) is a modular form of weight $4$ and level $3$ 
(resp., weight $6$ and level $2$) with Fourier coefficients in $\Q$, and 
$\gamma_2(z)$ and $\gamma_3(z)$ are modular functions of levels $3$ and $2$, respectively,
with Fourier coefficients in $\Q$.
Let $j(z)$ denote the usual $j$-function.
Weber (see for example p.\ 326 of \cite{Schertz}) showed
\begin{equation}
\label{g23a}
\gamma_2(z)^3 = j(z) \qquad \text{and} \qquad\gamma_3(z)^2 =j(z)-1728.
\end{equation}

If $F$ is a subfield of $\bar{\Q}$ or is a local field, let $\O_F$ denote its ring of integers.

If $F \subset \C$ is a number field, 
let $\id{F}$ denote its idele group, 
and let $F^{\ab}$ denote the maximal abelian extension of $F$ in $\C$.  
If $s \in \id{F}$ let $[s,F] \in \Gal(F^\ab/F)$ 
denote its global Artin symbol.  
If $w$ is a place of $F$ then $F_w$ will denote the completion of $F$ 
at $w$, and if $s \in \id{F}$ then $s_w \in F_w^\times$ will 
denote the $w$-component of $s$.  

By a prime of a number field $F$ we mean a prime ideal of $\O_F$.
If $\P$ is a prime of $F$, 
let $F^{\ab,\P}$ denote the maximal extension of 
$F$ in $F^{\ab}$ that is unramified at $\P$,
and if $a\in F^\times$, let 
$\ord_\P(a)$ be the power of $\P$ in the prime factorization of the fractional 
ideal $a\O_F$.
The Frobenius automorphism $\Fr_\P$
associated to $\P$ is the unique $\sigma\in\Gal(F^{\ab,\P}/F)$
such that 
$\sigma(x) \equiv x^{ \N_{F/\Q}(\P)} \pmod{\P\O_{F^{\ab,\P}}}$ for all $x\in\O_{F^{\ab,\P}}$.

Let $\R^+$ denote the multiplicative group of positive real numbers,
let $\GL_2^+(\R)$ (respectively, $\GL_2^+(\Q)$)
denote the subgroup of $\GL_2(\R)$  (respectively, $\GL_2(\Q)$) of 
elements with positive determinant,
 and let
$\GL_2^+(\ad{\Q})$ denote the subgroup of $\GL_2(\ad{\Q})$ consisting of 
elements whose $\infty$-com\-po\-nent has positive determinant.
Let
$$
\U = \GL_2^+(\R)  \times \prod_\ell \GL_2(\Zl) \subset \GL_2^+(\ad{\Q}).
$$
Recall that $g= \abcd \in \GL_2^+(\Q)$ acts on $\H$ by $g(z) = \frac{\alpha z+\beta}{\gamma z+\delta}$.

\begin{defn}
\label{Shimuraaction}
Shimura 
(see \cite{Shimura1978} or \S A5 of \cite{val}; see also \S 6.6 of \cite{Shimurabook} or \S1 of \cite{Rumely})
defined an action of $\GL_2^+(\ad{\Q})$ on the space of 
modular forms $f$ of weight $k$ with Fourier coefficients in $\Qab$, 
for every $k \in \Z$,
characterized by:
\begin{enumerate}
\item
the subgroup of $\GL_2^+(\ad{\Q})$ fixing $f$ is open,
\item
\label{fgz}
$f^g(z) = (\gamma z+\delta)^{-k}f(g(z))$ for every 
$g= \abcd \in \GL_2^+(\Q)$,
and
\item
\label{iotas}
if $s\in \R^+ \times \prod_\ell \Z_\ell^\times$ and 
$\iota(s):=\tinymatrix{1}{0}{0}{s^{-1}}$,
then $f^{\iota(s)}=f^{[s,\Q]}$, where
$[s,\Q]$ acts on $f$ by acting on the Fourier coefficients.
\end{enumerate}
\end{defn}

If $K$ is an imaginary quadratic field and
$\zzero \in K\cap\H$, let $q_{\zzero} : K \to \M_2(\Q)$ be the map defined by 
$$
q_{\zzero}(\mu)
\begin{pmatrix}\zzero \\ 1\end{pmatrix}
=
\begin{pmatrix}\mu \zzero \\ \mu\end{pmatrix}.
$$
Then $q_{\zzero}(K^\times)\subseteq \GL_2(\Q)$. 
Extend $q_{\zzero}$ to a map $q_{\zzero} : \ad{K} \to \M_2(\ad{\Q})$. 
Note that  for all $\mu\in\id{K}$, 
\begin{equation}
\label{qzdet}
\det(q_{\zzero}(\mu))=\N_{K/\Q}(\mu)
\end{equation}
so in particular $\det(q_{\zzero}(\mu)_\infty) = \mu_\infty \bar\mu_\infty > 0$, 
and therefore $q_{\zzero}(\id{K}) \subseteq \GL_2^+(\ad{\Q})$.

The following theorem is Theorem 6.31(i) of \cite{Shimurabook}. 

\begin{thm}[Shimura Reciprocity]
\label{sr}
Suppose $f$ is a modular function with Fourier coefficients in $\Qab$, 
$K$ is an imaginary quadratic field, 
$\zzero  \in K \cap \H$, and $f$ is defined and finite at $\zzero$.  
Then $f(\zzero) \in \Kab$, and if $s \in \id{K}$ then 
$$
f(\zzero)^{[s,K]} = f^{q_{\zzero}(s)^{-1}}(\zzero).
$$
\end{thm}

Let $\bmu_n:=\{ z\in\C : z^n=1\}$.

\begin{defn}
\label{legendrequartic}
Suppose $F \subset \C$ is a number field, $\bmu_n \subset F$, 
$\P$ is a prime of $F$ not dividing $n$, 
and $a \in F^\times$ is such that $n | \ord_\P(a)$.  
Then $F(a^{1/n}) \subset F^{\ab,\P}$ and we
define the $n$-th power symbol 
$$
\leg{a}{\P}_{n,F} := (a^{1/n})^{(\Fr_\P-1)} \in \bmu_n. 
$$
\end{defn}
Note that if $m|n$ then $\leg{a}{\P}_{m,F} = \leg{a^{n/m}}{\P}_{n,F}$.
If further $a\in\O_F - \P$, then 
$\leg{a}{\P}_{n,F} \in \bmu_n$ is characterized 
by the congruence
$$
\leg{a}{\P}_{n,F} \equiv a^{(\N_{F/\Q}(\P)-1)/n} \pmod{\P}.
$$
When $n=2$ this is the quadratic residue symbol, and it is $1$ if $a$ is
a square in $(\O_F/\P)^\times$ and $-1$  if $a$ is
a nonsquare in $(\O_F/\P)^\times$.

If $E : y^2=x^3+ax+b$ is an elliptic curve, its discriminant   
$\Delta(E)$ is $-16(4a^3+27b^2)$. 
By $\End(E)$ we mean
endomorphisms defined over an algebraic closure of the ground field.
When $E$ is an elliptic curve over $\C$, let $E[N]=\{ P\in E(\C) : NP=O\}$.

\begin{defn}
\label{esubdef}
When $\alpha\in\C^\times$ and $z\in\H$, define an elliptic curve over $\C$:
$$
\Et{\alpha}{z} : y^2 = x^3 - \alpha^2 \frac{\gamma_2(z)}{48}x 
    + \alpha^3 \frac{\gamma_3(z)}{864}.
$$
\end{defn}
Then:  
\begin{equation}
\label{jEalph}
j(\Et{\alpha}{z}) = j(z), \;
    \Delta(\Et{\alpha}{z}) = \alpha^6, \; \text{and} \; 
    \End_\C(\Et{\alpha}{z}) = \{\lambda \in \C : \lambda L_z \subseteq L_z\} .
\end{equation}
When $\alpha = 1$ we will often write simply $E_{z}$ instead of $\Et{1}{z}$.

If $K$ is an imaginary quadratic field, $\O$ is an order in $K$, and 
$\ell$ is a rational prime, let $\O_{\ell} := \O \otimes_\Z \Z_\ell$.
If $s \in \id{K}$, let 
$s_\ell$ denote the projection of $s$ in 
$(K \otimes_\Q \Q_\ell)^\times \subset \id{K}$.  

\begin{defn}
\label{Ogendef}
Suppose $K$ is an imaginary quadratic field,
$\O$ is an order in $K$, $F$ is a finite extension of $K$, and $\P$ is a prime of $F$.  
Let 
$$
V_\P = \{x \in F_\P^\times : \ord_\P(x) = 1\} \subset F_\P^\times \subset \id{F}.
$$
We define an {\em $(\O,F)$-good generator of $N_{F/K}(\P)$} 
to be an element $\lambda \in K^\times$ such that 
\begin{equation}
\label{OFgoodeqn}
\lambda^{-1}\N_{F/K}(V_\P) \subset K_\infty^\times \prod_\ell \O_{\ell}^\times.
\end{equation}
\end{defn}

\begin{lem}
\label{Ogenlem}
Let $K$, $\O$, $F$, and $\P$ be as in Definition \ref{Ogendef}.
If $\lambda$ is an $(\O,F)$-good generator of $N_{F/K}(\P)$, then:
\begin{enumerate}
\item 
$\lambda \in \O_K$ and $\lambda\O_K = \N_{F/K}(\P)$,
\item
if  $u \in \O_K^\times$, 
then $u\lambda$ is $(\O,F)$-good if and only if $u \in \O^\times$,
\item
if $\P \nmid 2$, then 
$\lambda \in \O_2^\times$.
\end{enumerate}
\end{lem}

\begin{proof}
Let $\p$ be the prime of $K$ below $\P$.  
Suppose $\q$ is a prime of $K$ and $t \in V_\P$.
By \eqref{OFgoodeqn},
   $\ord_\q(\lambda) = 0 = \ord_\q(\N_{F/K}(\P))$ if $\q \neq \p$, and
   $\ord_\p(\lambda) = \ord_\p(\N_{F/K}(t)) = \ord_\p(\N_{F/K}(\P))$, so
$\lambda \O_K = \N_{F/K}(\P)$, giving (i).
If $u\in\O^\times$, then clearly $u\lambda$ is $(\O,F)$-good. Conversely,
if $\lambda$ and $\lambda'$ are both $(\O,F)$-good generators of $N_{F/K}(\P)$,
 then their ratio is in $\O_\ell^\times$ for every $\ell$, so it is in $\O^\times$.
 This gives (ii).
Assume $\P \nmid 2$. Then $N_{F/K}(V_\P) \in\O_2^\times$. Thus by \eqref{OFgoodeqn}, 
$\lambda\in\O_2^\times$, giving (iii).
\end{proof}

\begin{rem}
In general, an $(\O,F)$-good generator of $N_{F/K}(\P)$ may not exist.  
We will show in Corollary \ref{EoCor} below that if there is an elliptic curve 
$E$ defined over $F$ with CM by $\O$ and with good reduction at $\P$, then 
$N_{F/K}(\P)$ has an $(\O,F)$-good generator, and if further 
$\P$ does not divide the conductor of the order $\O$, then 
$N_{F/K}(\P)$ has an $(\O,F)$-good generator in $\O$, and a generator
of $N_{F/K}(\P)$ is $(\O,F)$-good if and only if it is in $\O$.
By Lemma \ref{Ogenlem}(ii), if $K$ is not $\Q(i)$ or $\Q(\sqrt{-3})$ and 
there is an $(\O,F)$-good generator of $N_{F/K}(\P)$, then every generator 
of the ideal $\N_{F/K}(\P)$ is $(\O,F)$-good.
\end{rem}

\section{Some background results}
\label{CtHc}

In this section we state or work out the 
properties of the Weber functions and Dedekind's $\eta$-function that 
we need to compute Hecke characters.

Fix an imaginary quadratic field $K$ and fix 
$\tau \in \H \cap K$.  Let $\O_\tau$ be the order associated to the 
lattice $L_\tau = \Z + \Z \tau$, 
i.e., 
$$
\O_\tau = \{\alpha \in K : \alpha L_\tau \subseteq L_\tau\}. 
$$
The ring class field $H_\tau$ of $\O_\tau$ is the abelian extension 
of $K$ corresponding under class field theory to the subgroup 
$K^\times K_\infty^\times \prod_\ell \O_{\tau,\ell}^\times$
of $\id{K}$.
Then $H_\tau = K(j(\tau))$ (see p.~23 of \cite{Deuringbk} or  
Theorem 5.7 of \cite{Shimurabook}). 
If $\lambda \in \Ol{\ell}^\times \subset \id{K}$ then 
$q_\tau(\lambda) \in \GL_2(\Z_\ell) \subset \GL_2^+(\ad{\Q})$, and if 
$s \in K_\infty^\times \prod_\ell \Ol{\ell}^\times$
then $q_\tau(s) \in \U$.
Note that  
$s \in K_\infty^\times \prod_\ell \Ol{\ell}^\times$ 
if and only if $s_\ell \in \Ol{\ell}^\times$ for every $\ell$.

\begin{defn}
\label{phidef}
Let $\phi : \SL_2(\Z/4\Z) \to \bmu_4$ be the unique homomorphism that 
sends $\tinymatrix{1}{1}{0}{1}$ to $i$.  We will also view $\phi$ as 
a homomorphism $\SL_2(\Z_2) \to \bmu_4$ by composing with reduction modulo $4$.
We define a function $\eone : \Ol{2}^\times \to \bmu_4$ as follows.  
If $\lambda \in \Ol{2}^\times$ then 
$\tinymatrix{1}{0}{0}{\N_{K/\Q}(\lambda)^{-1}}q_\tau(\lambda) \in \SL_2(\Z_2)$
by \eqref{qzdet}, and we let   
$$
\eone(\lambda) = \phi\bigl(\tinymatrix{1}{0}{0}{\N_{K/\Q}(\lambda)^{-1}}q_\tau(\lambda)\bigr) \in \bmu_4.
$$
Then $\eone(\lambda)$ depends only on the reduction of $\lambda$ modulo $4 \Ol{2}$, 
so we will also view $\eone$ as a function from
$(\O_{\tau,2}/4\O_{\tau,2})^\times$ to $\bmu_4$.
Note that $(\O_{\tau,2}/4\O_{\tau,2})^\times = (\O_\tau/4\O_\tau)^\times$.
\end{defn}

\begin{lem}
\label{lemeta}
Suppose $s \in \id{K}$ is such that $s_\ell \in \Ol{\ell}^\times$ 
for every rational prime $\ell$.  Then
$$
(\eta^{q_{\tau}(s)})^6 = \eone(s_2) \eta^6.
$$
\end{lem}

\begin{proof}
Let $\rho = \eta^6$, and for every $g = \abcd \in \SL_2(\Z)$ define 
$\rho|_g$ by
$$
(\rho|_g)(z) = (\gamma z+\delta)^{-3}\rho(g(z)).
$$
Then $\rho$ is a modular form of weight $3$ and 
level $4$ with Fourier coefficients in $\Q$, and 
$\rho|_g = \phi(g)\rho$ for every $g \in \SL_2(\Z)$ 
(see for example \S1 of \cite{hrv}).

Let
$$
\U_4 = \{ v \in \U : v_2 - \tinymatrix{1}{0}{0}{1} \in 4 M_2(\Z_2)\}, 
\quad\text{and}\quad
w = \tinymatrix{1}{0}{0}{\N_{K/\Q}(s)^{-1}}.
$$
By \eqref{qzdet}, $w \cdot q_{\tau}(s) \in \U \cap \SL_2(\ad{\Q})$, 
so by Lemma 1.38 of \cite{Shimurabook} we can write
\begin{equation}
\label{0part}
w \cdot q_{\tau}(s) = v \cdot h
\end{equation}
with $v \in \U_4$ and $h \in \SL_2(\Z)$.  Since the 
Fourier coefficients of $\rho$ lie in $\Q$, 
Definition \ref{Shimuraaction}\ref{iotas} shows that $\rho^w = \rho$.  
Since $\rho$ has level $4$, Proposition 1.4 of \cite{Shimura1978} shows that 
$\rho^v = \rho$.  Thus (using Definition \ref{Shimuraaction}(ii))
\begin{equation}
\label{1part}
\rho^{q_{\tau}(s)}
    = \rho^{w \cdot q_{\tau}(s)}
    = \rho^{v \cdot h}
    = \rho^h
    =  \rho|_h = \phi(h) \rho.
\end{equation}
Since $\phi(h) = \phi(w_2 q_\tau(s_2)) = \eone(s_2)$, this proves the lemma.
\end{proof}

The next result is an application of Shimura's Reciprocity Law.
Its proof  is similar to Rumely's proof of part of Theorem 1 of \cite{Rumely}.

\begin{prop}
\label{ptrans}
Suppose $N \in \Z^+$, $F$ is a finite extension of $K$, 
$\P$ is a prime of $F$ not dividing $2N$,   
and $u \in N^{-1}\O_\tau/\O_ \tau$.  Then:
\begin{enumerate}
\itemsep=5pt
\item
\label{inHab}  
${\wp'(u;\tau)}/({(2 \pi i)^3 \eta(\tau)^6}) \in F^{\ab,\P}$,
\item
\label{frpart}
If $\lambda$ is an $(\O_\tau,F)$-good generator of $\N_{F/K}(\P)$, 
then
$$
\left(\frac{\wp'(u;\tau)}{(2 \pi i)^3 \eta(\tau)^6}\right)^{\Fr_\P} = 
    \;\eone(\lambda)^{-1}\frac{\wp'(\lambda u;\tau)}{(2 \pi i)^3 \eta(\tau)^6}.
$$
\end{enumerate}
\end{prop}

\begin{proof}
For $T\in\U$, let $T_N$ denote the image of $T$ in $\GL_2(\Z/N\Z)$.
If $(a,b) \in (N^{-1}\Z/\Z)^2$ (viewed as a row vector), define 
$$
f_{(a,b)}(z) = \frac{\wp'(az+b;z)}{(2\pi i)^3}.
$$
Then (see \S6.1 and \S6.2 of \cite{Shimurabook}, or p.\ 392 of \cite{Rumely}), 
\begin{enumerate}
\renewcommand{\theenumi}{(\alph{enumi})}
\item
$f_{(a,b)}$ is a modular form of weight $3$ with Fourier coefficients 
in $\Q^\ab$, 
\item
if $T \in \U$ then 
$(f_{(a,b)})^T = f_{(a,b)T_N}$.
\end{enumerate}
Let $\p$ be the prime of $K$ below $\P$, let 
$p$ be the prime of $\Q$ below $\P$, and write 
$u = a\tau + b$ with $a,b \in N^{-1}\Z/\Z$.  Then $f_{(a,b)}/\eta^6$ 
is a modular function with Fourier coefficients in $\Q^{\ab}$, 
and
$$
\frac{\wp'(u;\tau)}{(2 \pi i)^3 \eta(\tau)^6} = \frac{f_{(a,b)}(\tau)}{\eta(\tau)^6}.
$$
Suppose $t \in F_\P^\times$ and $\ord_\P(t) = 1$.  View $t \in \id{F}$,  
and let $s = \lambda^{-1}\N_{F/K}(t) \in \id{K}$.
Since $\lambda$ is an $(\O_\tau,F)$-good generator of $\N_{F/K}(\P)$, 
we have $s \in K_\infty^\times \prod_\ell \Ol{\ell}^\times$, so $q_{\tau}(s) \in \U$.

By Theorem \ref{sr}, 
${\wp'(u;\tau)}/({(2 \pi i)^3 \eta(\tau)^6}) \in K^{\ab}$ and 
\begin{equation}
\label{l1}
\left(\frac{\wp'(u;\tau)}{(2 \pi i)^3 \eta(\tau)^6}\right)^{[s,K]} 
    = \left(\frac{f_{(a,b)}(\tau)}{\eta(\tau)^6}\right)^{[s,K]}  
    = \frac{(f_{(a,b)})^{q_{\tau}(s)^{-1}}(\tau)}{(\eta^{q_{\tau}(s)^{-1}}(\tau))^6}. 
\end{equation}
Let $(a',b') := (a,b)q_{\tau}(s)^{-1}_N \in (N^{-1}\Z/\Z)^2$.  
Since $\P \nmid N$, we have $s_\ell = \lambda^{-1}$ for all $\ell|N$,
and so $q_{\tau}(s)_N^{-1} = q_{\tau}(\lambda)_N$. Thus in $\C/\O_\tau$,
\begin{multline*}
a'\tau + b' = (a',b') \begin{pmatrix}\tau \\ 1\end{pmatrix}
    = (a,b)q_{\tau}(s)_N^{-1} \begin{pmatrix}\tau \\ 1\end{pmatrix}\\
    = (a,b)q_{\tau}(\lambda) \begin{pmatrix}\tau \\ 1\end{pmatrix}
    = (a,b) \begin{pmatrix}\lambda\tau \\ \lambda\end{pmatrix}
    = \lambda u.
\end{multline*}
Using this and (b) above, 
\begin{equation}
\label{l2}
(f_{(a,b)})^{q_{\tau}(s)^{-1}}(\tau) = f_{(a,b)q_{\tau}(s)_N^{-1}}(\tau) 
    = f_{(a',b')}(\tau) = \frac{\wp'(\lambda u;\tau)}{(2 \pi i)^3}.
\end{equation}
Since $\P \nmid 2$, we have $s_2 = \lambda^{-1}$, so by 
Lemma \ref{lemeta},  
$$
(\eta^{q_{\tau}(s)^{-1}}(\tau))^6 = \eone(\lambda)\eta(\tau)^6.
$$
Combining this with \eqref{l1} and \eqref{l2} immediately gives
\begin{equation}
\label{pprimeeta}
\left(\frac{\wp'(u;\tau)}{(2 \pi i)^3 \eta(\tau)^6}\right)^{[s,K]} 
    = \eone(\lambda)^{-1}\frac{\wp'(\lambda u;\tau)}{(2 \pi i)^3 \eta(\tau)^6}.
\end{equation}
Since the right-hand side is independent of $t$ (recall that $s$ was 
defined in terms of $t \in F_\P^\times$), for every $r \in \O_{F,\P}^\times$ 
we have 
$$
\frac{\wp'(u;\tau)}{(2 \pi i)^3 \eta(\tau)^6}
    = \left(\frac{\wp'(u;\tau)}{(2 \pi i)^3 \eta(\tau)^6}\right)^{[\N_{F/K}(r),K]} 
    = \left(\frac{\wp'(u;\tau)}{(2 \pi i)^3 \eta(\tau)^6}\right)^{[r,F]}.
$$
Since 
$\{[r,F] : r \in \O_{F,\P}^\times\}$ is the inertia group at $\P$ 
in $\Gal(F^\ab/F)$, it follows that 
${\wp'(u;\tau)}/({(2 \pi i)^3 \eta(\tau)^6}) \in F^{\ab,\P}$,
giving  \ref{inHab}.
Let $L=K^{\ab}\cap F^{\ab,\P}$.
By class field theory, 
$$
\bigl.[s,K]\bigr|_{L} 
    = \bigl.[\N_{F/K}(t),K] \bigr|_{L} 
    = \bigl.[t,F] \bigr|_{L} 
    = \bigl.\Fr_\P \bigr|_{L}.
$$
This and \eqref{pprimeeta} give \ref{frpart}.
\end{proof}

\begin{lem}
\label{notjlem}
Let $D \in \Z_{<0}$ denote the discriminant of the order $\O_\tau$. Then:
\begin{enumerate}
\item
\label{g3sq}
$\gamma_2(\tau)^3, \gamma_3(\tau)^2 \in \Q({\jd(\tau)}) \subset H_\tau$;
\item
if $D$ is odd then $\sqrt{D}{\gamma_3(\tau)} \in \Q({j(\tau)}) \subset H_\tau$ 
and ${\gamma_3(\tau)} \in H_\tau$;
\item
if $D \equiv 4$ or $8 \pmod{16}$ then 
$\sqrt{-D}{\gamma_3(\tau)} \in \Q({j(\tau)}) \subset H_\tau$ and 
$i{\gamma_3(\tau)} \in H_\tau$;
\item
if $D \equiv 0$ or $12 \pmod{16}$ then $i \in H_\tau$.
\end{enumerate}
\end{lem}

\begin{proof}
Part  \ref{g3sq} follows from \eqref{g23a}. 
Let $\omega = (3+\sqrt{D})/2$ if $D$ is odd, and $\omega = \sqrt{D}/2$ 
if $D$ is even.  Then $L_\omega = \O_\tau$, so 
$\gamma_3(\omega)^2 = j(\omega) - 1728$ and $\gamma_3(\tau)^2 = j(\tau) - 1728$ 
are $\Gal(H_\tau/K)$-conjugates by Theorem 5.7 of \cite{Shimurabook}.  
Therefore it suffices to prove (ii), (iii), (iv) when $\tau$ is replaced by $\omega$.  In this case 
all three statements 
(except $D=-8$, which is easy to check)
are proved by Birch in \S6 of \cite{Birch} (who in turn 
says that they were either proved or noticed by Weber in \S\S125, 126, 134 of 
\cite{Weber}).
\end{proof}

\section{Computing the Hecke character}
\label{heckesect}

As before, fix an imaginary quadratic field $K$ and 
fix $\tau \in \H \cap K$. 
Theorem \ref{ptranscor} below is the key to our main results in 
\S\ref{Mt}.  For example, when $F = H_\tau$ 
it allows us to compute the Hecke character of $\Et{\alpha}{\tau}$ 
over $H_\tau$ whenever $\Et{\alpha}{\tau}$ is defined over $H_\tau$, even if 
$\alpha \notin H_\tau$ (i.e., even if $E_{\tau}$ is not defined over $H_\tau$).
We first state the basic properties we will need of the Hecke character.

\begin{prop}
\label{Eorder}
Suppose $E$ is an elliptic curve over a number field $F \supseteq K$, 
and $\O:=\End(E)$ is an order in $K$.
Let $B$ be the set of primes of $F$ where $E$ has bad reduction, 
and let $I(B)$ be the group of fractional ideals of $F$ supported 
outside of $B$.
Then there is a unique character $\psi = \psi_{E/F} : I(B) \to K^\times$,
called the {\em Hecke character} of $E$ over $F$, 
such that for every prime $\P$ of $F$ where $E$ has good reduction:
\begin{enumerate}
\item
\label{Eorderok}
$\psi(\P)\in\O_K$, and $\psi(\P)$ is an $(\O,F)$-good generator of $\N_{F/K}(\P)$;
\item
if  $\O=\Z+cp^r\O_K$, where $p$ is the residue characteristic of $\P$ 
and $p \nmid c$,  then 
$\psi(\P) \in \Z+c\O_K$;
\item
if $\P$ does not divide the conductor of $\O$ then $\psi(\P) \in \O$;  
\item
\label{Eorderpts}
$
|E(\O_F/\P)| = \N_{F/\Q}(\P) + 1 - \Tr_{K/\Q}(\psi(\P)).
$
\end{enumerate}
\end{prop}

\begin{proof}
Let $\psi_{\ad{}} : \id{F} \to \C^\times$ denote the Hecke character of $E$ over $F$ 
on ideles, as defined in \S 7.8 of \cite{Shimurabook}.  By Theorem 7.42 of 
\cite{Shimurabook}, $\psi_{\ad{}}$ is unramified at $\P$.
Then 
$\psi(\P) = \psi_{\ad{}}(t)$ where $t \in F_\P^\times \subset \id{F}$ is any element 
satisfying $\ord_\P(t) = 1$.  It follows from Proposition 7.40(ii) of \cite{Shimurabook} 
that $\psi(\P)/N_{F/K}(t) \in K_\infty^\times \prod_\ell \O_\ell^\times$, 
so $\psi(\P)$ is an $(\O,F)$-good generator of $N_{F/K}(\P)$
(in the sense of Definition \ref{Ogendef}). By Lemma \ref{Ogenlem}(i),
we have $\psi(\P)\in\O_K$,  giving (i).

For (ii), we follow a standard method as in, for example, the proof of
Theorem 12 in Chapter 13 of  \cite{lang}.
Let $\Eb$ denote the reduction of $E$ modulo $\P$, and let $p$ be the 
rational prime below $\P$.  It is shown in the proof 
of Theorem 7.42 of \cite{Shimurabook} that the image of $\psi(\P)$ under 
$$
K = \O \otimes \Q = \End(E) \otimes \Q \hookrightarrow \End(\Eb) \otimes \Q
$$
is the Frobenius endomorphism $\varphi \in \End(\Eb) \subset \End(\Eb) \otimes \Q$.  
Thus for every rational prime $\ell \ne p$, if $T_\ell$ denotes the $\ell$-adic Tate module 
we have a commutative diagram
$$
\xymatrix@C=40pt{
T_\ell(E)\otimes \Q \ar^{\psi(\P)}[r] \ar_{\cong}[d] & T_\ell(E)\otimes \Q \ar^{\cong}[d] \\
T_\ell(\Eb)\otimes \Q \ar^{\varphi}[r] & T_\ell(\Eb)\otimes \Q
}
$$
where the vertical maps are induced by the reduction isomorphism 
$T_\ell(E) \isom T_\ell(\Eb)$.
Since $\varphi \in \End(\Eb)$, we have $\varphi(T_\ell(\Eb)) \subseteq T_\ell(\Eb)$.  
Thus by Theorem 5 of \cite{serretate}, $\psi(\P) \in \O_\ell$  
for all $\ell \ne p$.  Thus
$$
\psi(\P) \in \O_K \bigcap_{\ell \ne p} \O_\ell = \Z+c\O_K.
$$
This gives (ii). 
If $\P$ does not divide the conductor $cp^r$ of $\O$ (i.e., $r = 0$), then 
$\Z+c\O_K = \Z+cp^r\O_K = \O$, giving (iii).

For \ref{Eorderpts}, see for example Corollary II.10.4.1 of  \cite{Silverman} for
the case where $\O$ is the maximal order $\O_K$, and see
Theorem 7.42 of \cite{Shimurabook} for the general case.
\end{proof}

\begin{cor}
\label{EoCor}
Suppose that $F$ is a number field containing $K$, $\P$ is a prime of $F$, 
and $\O$ is an order in $K$.  
If there is an elliptic curve $E$ defined over $F$ with CM by $\O$ and 
with good reduction at $\P$, then:
\begin{enumerate}
\item
 $N_{F/K}(\P)$ has an $(\O,F)$-good generator;
\item
if $\P$ does not divide the 
conductor of the order $\O$, then:
\begin{enumerate}
\item
$N_{F/K}(\P)$ has a generator in $\O$,
\item
 a generator of 
$N_{F/K}(\P)$ is $(\O,F)$-good if and only if it lies in $\O$.
\end{enumerate} 
\end{enumerate} 
\end{cor}

\begin{proof}
By Proposition \ref{Eorder}(i), $\psi(\P)$ is an $(\O,F)$-good generator of $N_{F/K}(\P)$, 
where $\psi$ is the Hecke character of $E$.
If $\P$ does not divide the conductor of $\O$, then
$\psi(\P)\in\O$ by Proposition \ref{Eorder}(iii).
Part (b) now follows from Lemma \ref{Ogenlem}(ii).
\end{proof}

Next we give an example in which $N_{F/K}(\P)$ has no generators in $\O$,
under the hypotheses in Corollary \ref{EoCor} (and Theorem \ref{main}),
so $\psi(\P)\notin\O$.
This is why we take an  $(\O,F)$-good generator of $N_{F/K}(\P)$, which
always exists by Corollary \ref{EoCor}(i), rather than a generator in $\O$.

\begin{exa}
\label{psimythex}
Let $K = \Q(\sqrt{-11})$. Then $\O_K=\Z[\beta]$ where $\beta=(1+\sqrt{-11})/2\in\O_K$. 
Let $\O = \Z + 3\O_K$,  the order of conductor $3$ in $\O_K$.
Then $3 = \beta\bar{\beta}$,
$$
j(\O) = j(3\beta) = -18808030478336 - 3274057859072\sqrt{33},
$$ 
and
$H_\O = K(j(\O)) = K(\sqrt{33}) = K(\sqrt{-3})$.
Let $E$ be the elliptic curve 
$$y^2 + y = x^3 -\frac{(7 + \sqrt{33})}{2}x^2 - \frac{(2487 + 433\sqrt{33})}{2}x - 21416 - 3728\sqrt{33}.$$
Then $E$ is defined over $F := H_\O$.
Since $j(E) = j(\O)$, $E$ has CM by $\O$.
The discriminant of $E$ is the unit $-23 - 4\sqrt{33}$, 
so $E$ has good reduction everywhere.
Let $\P$ be a prime of $F$ above $\beta$.
Since $\P$  is totally ramified in the extension $F/K$,
we have $N_{F/K}(\P) = \beta\O_K$, which has no generators in $\O$.
Therefore, $\psi(\P)\notin\O$.
Note that the reduction of $E$ mod $\P$ has CM by $\O_K$.
\end{exa}

Recall $\eone$ from Definition \ref{phidef}.

\begin{thm}
\label{ptranscor}
Suppose $K$ is an imaginary quadratic field, 
$\tau \in \H \cap K$, and $\O_\tau^\times = \{\pm1\}$.
Suppose $F$ is a number field containing $K$, and
$\alpha\in\C^\times$ is such that $\alpha^2 {\gammatwod}(\tau), \alpha^3 {\gammathreed}(\tau)\in F$.
Let $\psi$ 
be the Hecke character of $\Et{\alpha}{\tau}$ over $F$.  If $\P$ is a prime ideal
of $F$ where $\Et{\alpha}{\tau}$ has good reduction, $\P \nmid 2$, 
and $\lambda$ is an $(\O_\tau,F)$-good generator of $\N_{F/K}(\P)$, then:
\begin{enumerate}
\item
\label{alpha6part}
$\alpha^6 \in F$,
\item
\label{4ordpart}
$4 \mid \ord_\P(\alpha^6)$, 
\item
$\psi(\P) = \pm \lambda$, and
\item
$\psi(\P) = \eone(\lambda) (\alpha^{9/2})^{(\Fr_\P-1)} \lambda 
= \eone(\lambda)^{-1} (\alpha^{3/2})^{(\Fr_\P-1)} \lambda$.
\end{enumerate}
\end{thm}

\begin{proof}
Let $j = j(\tau)$, $\gamma_2 = \gamma_2(\tau)$, 
and $\gamma_3 = \gamma_3(\tau)$. 
Note that $H_\tau = K(j) = K(j(\Et{\alpha}{\tau})) \subseteq F$.  
Since $\gammatwod^3$ and $\gammathreed^2 \in H_\tau$ (by Lemma \ref{notjlem}(i)), 
and $\gammatwod^3$ and $\gammathreed^2$ 
cannot both be zero (by \eqref{g23a}), we have \ref{alpha6part}.

By \eqref{jEalph}, $\End_{\C}(\Et{\alpha}{\tau}) = \O_\tau$.
The map $t : \C/L_\tau \to \Et{\alpha}{\tau}(\C)$ defined by 
$$
t(u) = (\alpha \wp(u;\tau)/((2 \pi i)^2 \eta(\tau)^4), 
    \alpha^{3/2}\wp'(u;\tau)/((2 \pi i)^3 \eta(\tau)^6))
$$
is an $\O_\tau$-module isomorphism.
Suppose  $N\in\Z^+$ is prime to $\P$
and suppose $u \in N^{-1}\O_ \tau/\O_ \tau = (\C/\O_ \tau)[N]$.  Then 
$t(u) \in \Et{\alpha}{\tau}[N]$.  
Since $\Et{\alpha}{\tau}$ has good reduction at $\P$
and $\P \nmid N$, the coordinates of $t(u)$ 
generate an extension of $F$ that is unramified at $\P$.  
By Proposition \ref{ptrans}\ref{inHab} it follows that $F(\alpha^{3/2})/F$ 
is unramified at $\P$, and since $(\alpha^{3/2})^4 = \alpha^6 \in F$
this proves  \ref{4ordpart}.

By Proposition 7.40(2) of \cite{Shimurabook}, 
$t(u)^{\Fr_\P} =  t(\psi(\P)u)$.
Taking $y$-coordinates and applying Proposition \ref{ptrans}(ii) gives
$$
\alpha^{3/2}\frac{\wp'(\psi(\P)u;\tau)}{(2 \pi i)^3 \eta(\tau)^6} = 
    \left(\alpha^{3/2}\frac{\wp'(u;\tau)}{(2 \pi i)^3 \eta(\tau)^6}\right)^{\Fr_\P} \\
    = (\alpha^{3/2})^{\Fr_\P}
            \eone(\lambda)^{-1}       
       \frac{\wp'(\lambda u;\tau)}{(2 \pi i)^3 \eta(\tau)^6}
$$
so
\begin{equation}
\label{lp}
\wp'(\psi(\P)u;\tau) = \wp'(\lambda u;\tau) 
    (\alpha^{3/2})^{(\Fr_\P-1)}
    \eone(\lambda)^{-1}.
\end{equation}
Since \eqref{lp} holds for a dense set of $u \in \C$, it holds for 
every $u \in \C$ by continuity.
The left side of \eqref{lp} has poles exactly at all $u\in \psi(\P)^{-1}L_\tau$
while the right side has poles exactly at all $u\in \lambda^{-1}L_\tau$.
Thus $\psi(\P)/\lambda \in \O_\tau^\times = \{\pm 1\}$, giving (iii).
Since $\wp'$ is an odd function, 
\begin{equation}
\label{lp2}
\wp'(\psi(\P)u;\tau) = \wp'((\psi(\P)/\lambda)\lambda u;\tau) 
    = (\psi(\P)/\lambda) \wp'(\lambda u;\tau)
\end{equation}
for all $u \in \C$. 
Comparing this with \eqref{lp} gives
\begin{equation}
\label{ptranscora}
\psi(\P)/\lambda = 
  \eone(\lambda)^{-1} (\alpha^{3/2})^{(\Fr_\P-1)}\in\{\pm1\} .
\end{equation}
Since $\alpha^6\in F$ by \ref{alpha6part},
we have $(\alpha^6)^{(\Fr_\P-1)}=1$ and thus
$$
\eone(\lambda)^{-1} (\alpha^{3/2})^{(\Fr_\P-1)}=
  \eone(\lambda) (\alpha^{-3/2})^{(\Fr_\P-1)}=
    \eone(\lambda) (\alpha^{9/2})^{(\Fr_\P-1)}.
$$
Combining this with \eqref{ptranscora} proves (iv).
\end{proof}

\section{Explicit formulas for Hecke characters and point counting}
\label{Mt}

The main results of this paper are Theorem \ref{main} and Corollary  \ref{exrem}.

If $K$ is an imaginary quadratic field and $\tau\in\H \cap K$,
let $D(\tau)$ denote
the discriminant of the order $\O_\tau$
(so $D(\tau)=B^2-4AC \equiv 0$ or $1\pmod{4}$ where $A\tau^2+B\tau+C=0$ with
$A, B, C \in\Z$ and $\gcd(A,B,C)=1$).

\begin{defn}
\label{e2def}
With $\tau$ as above and using  $\eone$ of Definition \ref{phidef},
define a map $\etwo : (\O_\tau/4\O_\tau)^\times \to \bmu_4$ by 
$$
\etwo(\lambda) = 
\begin{cases}
i^{ (\N_{K/\Q}(\lambda) - 1)/2}\eone(\lambda) & \text{if $D(\tau) \equiv 4$ or $8 \pmod{16}$}, \\
\eone(\lambda) & \text{otherwise}.
\end{cases}
$$
\end{defn}
We will give $\etwo$ in a concrete and explicit way in \S\ref{esect}.

Recall the quadratic and quartic symbols 
$\legtF{a}{\P}{F}$ and $\legfF{a}{\P}{F}$ of 
Definition \ref{legendrequartic}.  

\begin{rem}
In Theorem \ref{main} below, if $K$ is not $\Q(i)$ or $\Q(\sqrt{-3})$, then by 
Lemma \ref{Ogenlem}(ii) and 
Proposition \ref{Eorder}(i), every generator of the principal ideal
$N_{F/K}(\P)$ is $(\O,F)$-good.  Thus 
in this case
the hypothesis ``let $\lambda$ be an $(\O,F)$-good generator of $N_{F/K}(\P)$'' can be replaced 
by ``let $\lambda$ be a generator of $N_{F/K}(\P)$''.  For arbitrary $K$, if $\P$ 
does not divide the conductor of the order $\O$, then by Corollary \ref{EoCor}(ii), 
the hypothesis ``let $\lambda$ be an $(\O,F)$-good generator of $N_{F/K}(\P)$'' can be replaced 
by  ``let $\lambda$ be a generator of $N_{F/K}(\P)$ in $\O$''.
The same simplifications apply to Corollary \ref{exrem}.
\end{rem}

\begin{thm}
\label{main}
Suppose $E : y^2 = x^3 + ax + b$ is an elliptic curve over a number field $F$, 
and
$\O:=\End(E)$ is an order in an imaginary quadratic field $K \subseteq F$. 
Assume $\O^\times = \{\pm 1\}$.
Take any $\tau\in \H \cap K$ such that $j(E) = j(\tau)$.  
Suppose $\P$ is a prime of $F$, not dividing $2$, 
where $E$ has good reduction. Let $\lambda$ be an 
$(\O,F)$-good generator of $N_{F/K}(\P)$, let $q = \N_{F/\Q}(\P)$,
let $\psi$ denote the Hecke character of $E$ over $F$,
let $D$ be the discriminant of $\O$,
and  let $j=j(\tau)$, $\gamma_2 = \gamma_2(\tau)$, and $\gamma_3 = \gamma_3(\tau)$
(so $\O=\O_\tau$).   
Then:
\begin{enumerate}
\item
If $D$ is odd, then $\gammathreed \in F$, $\ord_\P(6b{\gammathreed})$ is even,
$$
\psi(\P) = \legtF{6b{\gammathreed}}{\P}{F}\etwo(\lambda)\lambda,
$$
 and 
$
|E(\O_F/\P)| = q + 1 - \legtF{6b{\gammathreed}}{\P}{F}\etwo(\lambda)\Tr_{K/\Q}(\lambda).
$
\item
If $D \equiv 4$ or $8 \pmod{16}$, then $i \gammathreed \in F$, 
$\ord_\P(-6bi{\gammathreed})$ is even, 
$$
\psi(\P) = \legtF{-6bi{\gammathreed}}{\P}{F}\etwo(\lambda)\lambda,
$$ 
and 
$
|E(\O_F/\P)| = q + 1 - \legtF{-6bi{\gammathreed}}{\P}{F}\etwo(\lambda)\Tr_{K/\Q}(\lambda).
$
\item
If $D \equiv 0$ or $12 \pmod{16}$, then $i \in F$, $4 \mid \ord_\P(6^2b^2({\jd}-1728))$, 
$$
\psi(\P) = \legfF{6^2b^2({\jd}-1728)}{\P}{F}\etwo(\lambda) \lambda,
$$ and 
$
|E(\O_F/\P)| = q + 1 - \legfF{6^2b^2({\jd}-1728)}{\P}{F}\etwo(\lambda) \Tr_{K/\Q}(\lambda).
$
\end{enumerate}
\end{thm}

\begin{proof}
The choice of $\tau$ implies that $\O=\O_\tau$.
Let 
$
\mm= 2^7 3^4 a^2b/(4a^3+27b^2)\in F^\times.
$
The map $(x,y)\mapsto (\mm^2 x,\mm^3 y)$
defines an isomorphism over $F$ from $E$ to the curve 
$y^2 = x^3+\mm^4ax+\mm^6b$.
The latter is 
$$
\textstyle
y^2 = x^3-\frac{3}{4}b^2j^3(j-1728)x+\frac{1}{4}b^3j^4(j-1728)^2,
$$
which is $\Et{\alpha}{\tau}$ with $\alpha := 6b\gamma_2^4\gamma_3$,
since 
$$
\gamma_2^3=j=j(E)=2^8 3^3 a^3/(4a^3+27b^2), \quad
\gamma_3^2 = j-1728=-2^63^6b^2/(4a^3+27b^2).
$$ 
Thus $E$ is isomorphic over $F$ to $\Et{\alpha}{\tau}$, so they have
the same Hecke character $\psi$ over $F$.  
Since $j \in F$, we have $H_\tau \subseteq F$.

{\em Case 1.} Suppose $D$ is odd.   Then 
$\alpha^9 = 6^9 b^9 \gamma_3j^{12}(j-1728)^4 \in F^\times$ 
by Lemma \ref{notjlem}(ii), 
and $\ord_\P(\alpha^9)$ is even by Theorem \ref{ptranscor}(ii), so
\begin{equation}
\label{a93}
(\alpha^{9/2})^{(\Fr_\P-1)} = \legtF{\alpha^9}{\P}{F} = \legtF{6b\gamma_3}{\P}{F}.
\end{equation}

{\em Case 2.} Suppose $D \equiv 4$ or $8 \pmod{16}$. Then 
$i\alpha^9 = 6^9 b^9 i\gamma_3j^{12}(j-1728)^4 \in F^\times$
by Lemma \ref{notjlem}(iii), 
and $\ord_\P(i\alpha^9)$ is even by Theorem \ref{ptranscor}(ii).  
If $\zeta \in \bmu_8$, then $\zeta^{(\Fr_\P-1)} = \zeta^{(q-1)}$.  Thus,
\begin{equation}
\label{a92}
(\alpha^{9/2})^{(\Fr_\P-1)} =
i^{(q-1)/2}\legtF{-i\alpha^9}{\P}{F} 
    = i^{(q-1)/2}\legtF{-6bi\gamma_3}{\P}{F}.
\end{equation}

{\em Case 3.} Suppose $D \equiv 0$ or $12 \pmod{16}$. Then
$$
\alpha^{18} = 6^{18} b^{18} j^{24}(j-1728)^9 \in F^\times,
$$
$i \in F$ by Lemma \ref{notjlem}(iv), 
and $4 \mid \ord_\P(\alpha^{18})$ by Theorem \ref{ptranscor}(ii).  
It follows that
\begin{equation}
\label{a91}
(\alpha^{9/2})^{(\Fr_\P-1)} = \legfF{\alpha^{18}}{\P}{F} = \legfF{6^2 b^2 (j-1728)}{\P}{F}.
\end{equation}

The desired formulas for $\psi(\P)$ now 
follow from 
Theorem \ref{ptranscor}(iv) along with \eqref{a93}, \eqref{a92}, \eqref{a91}, and
Definition \ref{e2def}. 
By Theorem \ref{ptranscor}(iii), 
$\psi(\P)/\lambda \in \{\pm1\}$, so 
$$\Tr_{K/\Q}(\psi(\P)) = (\psi(\P)/\lambda)\Tr_{K/\Q}(\lambda).$$ 
The desired formulas for $|E(\O_F/\P)|$
now follow from Proposition \ref{Eorder}\ref{Eorderpts}.
\end{proof}

\begin{cor}
\label{exrem}
Suppose $K$ is an imaginary quadratic field, 
$\tau \in \H \cap K$, and $\O_\tau^\times = \{\pm 1\}$.
Suppose $F$ is a finite extension of $H_\tau$ and $\beta \in F^\times$. 
With  $j:=j(\tau)$, $\gamma_2 := \gamma_2(\tau)$, and $\gamma_3 := \gamma_3(\tau)$,
let $E$ be the elliptic curve given by the following table, depending on 
$D(\tau) \pmod{16}:$
$$
\renewcommand{\arraystretch}{1.75}
\begin{array}{|l|rl|}
\hline
D(\tau)\qquad & E &\\
\hline
\text{odd} & \Et{\beta\gammatwod^4}{\tau} : &y^2 = x^3 -\frac{\beta^2 {\jd}^3}{48} x 
    +\frac{\beta^3 {\gammathreed}{\jd}^4}{864}\\
\hline \text{$4$ or $8 \pmod{16}$} & \Et{\beta i\gammatwod^4}{\tau} : 
    &y^2 = x^3 + \frac{\beta^2 {\jd}^3}{48}x 
        - \frac{\beta^3 i{\gammathreed}{\jd}^4}{864}\\ 
\hline \text{$0$ or $12 \pmod{16}$} & \Et{\beta\gammatwod^4{\gammathreed}}{\tau} : 
    &y^2 = x^3 - \frac{\beta^2 {\jd}^3({\jd}-1728)}{48}x 
        + \frac{\beta^3 {\jd}^4({\jd}-1728)^2}{864}\\
\hline
\end{array}
$$
Suppose $\P$ is a prime of $F$, not dividing $2$, where $E$ has good reduction.
Suppose $\lambda$ is an $(\O_\tau,F)$-good generator of $N_{F/K}(\P)$.
Let $q = \N_{F/\Q}(\P)$.
Then:
\begin{enumerate}
\item
$E$ is defined over $F$, $\End(E)=\O_\tau$, and $j(E) = {\jd}$;
\item
if $D(\tau)$ is odd or 
$D(\tau) \equiv 4$ or $8 \pmod{16}$, and
$\psi$ is the Hecke character of $E$ over $F$, then $\ord_\P(\beta)$ is even, 
$\psi(\P) = \legtF{\beta}{\P}{F}\etwo(\lambda)\lambda$, and
$$
|E(\O_F/\P)| = q+1-\legtF{\beta}{\P}{F}\etwo(\lambda) \Tr_{K/\Q}(\lambda);
$$
\item
if $D(\tau) \equiv 0$ or $12 \pmod{16}$, and
$\psi$ is the Hecke character of $E$ over $F$, 
then $4$ divides $\ord_\P(\beta^2({\jd}-1728))$, 
$
\psi(\P) = \legfF{\beta^2({\jd}-1728)}{\P}{F}\etwo(\lambda)\lambda,
$
and 
$$
|E(\O_F/\P)| = q+1- \legfF{\beta^2({\jd}-1728)}{\P}{F}\etwo(\lambda) \Tr_{K/\Q}(\lambda).
$$
\end{enumerate}
\end{cor}  

\begin{proof}
By  Lemma \ref{notjlem}, $E$ is defined over $F$.
By \eqref{jEalph}, $j(E) = {\jd}=j(\tau)$, so 
$\End(E)=\O_\tau$.
Now (ii) and (iii) follow directly from Theorem \ref{main}, using
the fact that $864 = 6\cdot 12^2$.
\end{proof}  

\begin{rem}
In Theorem \ref{main} we exclude the cases where $\O_\tau^\times$ is larger than 
$\{\pm1\}$.  This excludes precisely those $\tau$ with $j(\tau) = 1728$ 
(i.e., $\O_\tau = \Z[i]$; i.e., $D(\tau) = -4$) or
$j(\tau) = 0$ (i.e., $\O_\tau = \Z[e^{2\pi i/3}]$; i.e., $D(\tau) = -3$). 
For completeness we include these cases in the next two results, which follow easily from classical 
results that go back to Gauss (see for example p.~318 of \cite{cox}).

\begin{thm}
\label{D=-4thm}
Suppose $F$ is a number field containing $i$. Suppose $a\in F^\times$, and 
$E$ is the elliptic curve
$y^2 = x^3 - ax$.  Let $\psi$ denote the Hecke character of $E$ over $F$.   
Suppose $\P$ is a prime of $F$, not dividing $2$, 
where $E$ has good reduction. Let $\lambda \in \Z[i]$ be the 
generator of the principal ideal $N_{F/\Q(i)}(\P)$ 
congruent to $1 \pmod{2+2i}$, and let $q = \N_{F/\Q}(\P)$.  
Then $4 \mid \ord_\P(a)$,
$$
\psi(\P) = \legfF{a}{\P}{F}^{-1}\lambda, \quad \text{ and }  \quad
|E(\O_F/\P)| = q + 1 - \Tr_{K/\Q}(\legfF{a}{\P}{F}^{-1}\lambda).
$$
\end{thm}

\begin{thm}
\label{D=-3thm}
Suppose $F$ is a number field containing $\sqrt{-3}$. Suppose $b\in F^\times$, and 
$E$ is the elliptic curve
$y^2 = x^3 + 16b$.  Let $\psi$ denote the Hecke character of $E$ over $F$.   
Suppose $\P$ is a prime of $F$, not dividing $6$, 
where $E$ has good reduction. Let $\lambda \in \Z[e^{2\pi i/3}]$ be the  
generator of the principal ideal $N_{F/\Q(\sqrt{-3})}(\P)$  
congruent to $1 \pmod{3}$, and let $q = \N_{F/\Q}(\P)$.  
Then $6 \mid \ord_\P(b)$,
$$
\psi(\P) = \leg{b}{\P}_{\hskip-1pt 6,F}^{-1}\lambda,
\quad \text{ and }  \quad
|E(\O_F/\P)| = q + 1 - \Tr_{K/\Q}(\leg{b}{\P}_{\hskip-1pt 6,F}^{-1}\lambda).
$$
\end{thm}
\end{rem}

\section{Computing $\etwo$}
\label{esect}

In order to make Theorem \ref{main} and Corollary \ref{exrem} explicit, it is 
necessary to compute the function $\etwo$. 
For any given $\tau$, this is a simple computation,
following a method described  (for example) in \S1 of \cite{hrv} 
(see the proofs of Lemma \ref{lemphi} and Proposition \ref{eprop}  below).

Suppose $\O$ is an arbitrary order in an imaginary quadratic field $K$
and define $\td$ as in \eqref{zddef} below.
Proposition \ref{eprop} below gives the explicit values of the function $\etd$.
Suppose $E$ is an elliptic curve over $F \supseteq K$.
If  $j(E)=j(\O)$ ($=j(\td)$), then 
Theorem \ref{main} and Proposition \ref{eprop} together give
explicit formulas for the number of points on the reductions of $E$.
When $\O=\O_K$, this gives Theorem \ref{mainspec}.
Under the more general hypotheses in Theorem \ref{main}
(i.e., $j(E)=j(\a)$ for a proper $\O$-ideal $\a$), take any $\tau$ satisfying
the conclusion of Lemma \ref{trs}(i) below. Then Lemma \ref{trs}(ii) and
Proposition \ref{eprop} together give an explicit  value for the
$\etwo(\lambda)$ that occurs in Theorem \ref{main} and Corollary \ref{exrem}.

Throughout this section,
suppose $D$ is the discriminant of an order $\O$ in an imaginary quadratic field $K=\Q(\sqrt{D})$ 
(i.e., $D$ is a negative integer and $D \equiv 0$ or $1 \pmod{4}$). 
Define a positive integer $d$ by 
\begin{equation}
\label{ddef}
d = 
\begin{cases}
-D & \text{if $D$ is odd} \\
-D/4 & \text{if $D$ is even}
\end{cases}
\end{equation}
and let $\sqrt{-d}$ denote the square root of $-d$ in $\H$.  Then 
$K = \Q(\sqrt{-d})$, and we
define $\td \in \H \cap K$ by the following table:
\begin{equation}
\label{zddef}
\renewcommand{\arraystretch}{1.75}
\begin{array}{|r|c|c|c|c|}
\hline
D :  & 1\hskip-7pt\pmod{8} & 5\hskip-7pt\pmod{8} & \text{$4$ or $8\hskip-7pt\pmod{32}$} & \text{otherwise} \\
\hline
\td: & \frac{-3+\sqrt{-d}}{2} & \frac{3+\sqrt{-d}}{2} & 3 + \sqrt{-d} & \sqrt{-d} \\
\hline
\end{array}
\end{equation}
Then $\O = \O_{\td} = L_{\td} = \Z + \Z \td$ and $j(\O)=j(\td)$.  

The function $\etwo$ was defined in terms of the map $\phi$ of Definition \ref{phidef}.
A strategy for computing values of $\phi$ is given in  \S1 of \cite{hrv}.
We state the relevant ideas in the next lemma, and use them  below.

\begin{lem}
\label{lemphi}
Suppose $M \in \SL_2(\Z/4\Z)$ and $k\in\Z$. 
Let $C$ denote the commutator subgroup of $\SL_2(\Z/4\Z)$.
Then:
\begin{enumerate}
\item
$\phi(M)=i^k$ if and only if $\tinymatrix{1}{-k}{0}{1} M \in C$,
\item
$\phi\bigl(\tinymatrix{-1}{0}{0}{1} M \tinymatrix{-1}{0}{0}{1}^{-1}\bigr) 
    = \shortoverline{\phi(M)}$.
\end{enumerate}
\end{lem}

\begin{proof}
The explicit description of $C$ (see p.~498 of \cite{hrv})
shows that $\tinymatrix{1}{1}{0}{1}$ generates $\SL_2(\Z/4\Z)/C$.
Thus, given $M$, there is a unique $k\in\Z/4\Z$ so that $M_k := \tinymatrix{1}{-k}{0}{1} M \in C$.
Then $\phi(M_k) = 1$ (since $\bmu_4$ is abelian), so $\phi(M)=i^k$.  
Now (i) follows.
Part (ii) follows from (i) and the fact that 
$\tinymatrix{1}{-k}{0}{1} \tinymatrix{-1}{0}{0}{1} =\tinymatrix{-1}{0}{0}{1}\tinymatrix{1}{k}{0}{1}$.
\end{proof}
 
\begin{prop}
\label{eprop}
The map $\etd : \O_2^\times \to \bmu_4$ is given
by the following tables.

\medskip\noindent
If $D$ is odd: \nopagebreak

\medskip~
$\arraycolsep=3pt
\renewcommand{\arraystretch}{1.25}
\begin{array}{|r|c|c|}
\hline
\lambda^3~(\mathrm{mod}~4):  & 1, -\sqrt{-d} & -1, \sqrt{-d} \\
\hline
\etd(\lambda): & 1 & -1 \\
\hline
\end{array}$

\medskip\noindent
If $D \equiv 4 \pmod{16}:$ \nopagebreak

\medskip~
$\arraycolsep=3pt
\renewcommand{\arraystretch}{1.25}
\begin{array}{|r|c|c|}
\hline
\lambda~(\mathrm{mod}~4):  & 1, \sqrt{-d}, -1+2\sqrt{-d}, 2-\sqrt{-d} & 
    -1, -\sqrt{-d}, 1+2\sqrt{-d}, 2+\sqrt{-d} \\
\hline
\etd(\lambda): & 1 & -1 \\
\hline
\end{array}$

\medskip\noindent
If $D \equiv 8 \pmod{16}:$ \nopagebreak

\medskip~
$\renewcommand{\arraystretch}{1.25}
\begin{array}{|r|c|c|}
\hline
\lambda~(\mathrm{mod}~4):  & 1, -1+2\sqrt{-d}, \pm1+\sqrt{-d} 
    & -1, 1+2\sqrt{-d}, \pm1-\sqrt{-d} \\
\hline
\etd(\lambda): & 1 & -1  \\
\hline
\end{array}$

\medskip\noindent
If $D \equiv 12 \pmod{16}:$ \nopagebreak

\medskip~
$\arraycolsep=3pt
\renewcommand{\arraystretch}{1.25}
\begin{array}{|r|c|c|c|c|}
\hline
\lambda~(\mathrm{mod}~4):  & 1, 1+2\sqrt{-d} & 2+\sqrt{-d},\sqrt{-d} 
    & -1,-1+2\sqrt{-d} & 2-\sqrt{-d}, -\sqrt{-d} \\
\hline
\etd(\lambda): & 1 & i & -1 & -i \\
\hline
\end{array}$

\medskip\noindent
If $D \equiv 0 \pmod{16}:$ \nopagebreak

\medskip~
$\arraycolsep=3pt
\renewcommand{\arraystretch}{1.25}
\begin{array}{|r|c|c|c|c|}
\hline
\lambda~(\mathrm{mod}~4):  & 1,-1+2\sqrt{-d} & \pm1-\sqrt{-d} 
    & -1,1+2\sqrt{-d} & \pm1+\sqrt{-d} \\
\hline
\etd(\lambda): & 1 & i & -1 & -i \\
\hline
\end{array}$
\end{prop}

\begin{proof}
Since $\etwo$ is a simple modification of $\eone$ (Definition \ref{e2def}), it suffices 
to compute $\etdtilde(\lambda)$.  By Definition \ref{phidef},
\begin{equation}
\label{mathver}
\etdtilde(\lambda) = \phi\bigl(\tinymatrix{1}{0}{0}{\N_{K/\Q}(\lambda)^{-1}} q_{\td}(\lambda)\bigr).
\end{equation}
We follow the strategy for computing values of $\phi$ described in 
\S1 of \cite{hrv} (and Lemma \ref{lemphi} above).  
Find $k\in\{0,1,2,3\}$ such that 
$\tinymatrix{1}{-k}{0}{1} \tinymatrix{1}{0}{0}{\N_{K/\Q}(\lambda)^{-1}} q_{\td}(\lambda)$ is in  
the commutator subgroup of $\SL_2(\Z/4\Z)$
(given explicitly on p.\ 498 of \cite{hrv}).
Then $\etdtilde(\lambda)=i^k$ by Lemma \ref{lemphi}(i) and \eqref{mathver}.
We carried out this computation in Mathematica, and obtained the values in the tables.
\end{proof}

\begin{rem}
The discriminants of maximal orders in imaginary quadratic fields are
exactly the negative integers $D$ such that either $D$ is squarefree and
$D\equiv 1\pmod{4}$, or $D=-4d$ with $d\in\Z^+$ squarefree and $d\equiv 1$ or $2\pmod{4}$.
So if $D$ is the discriminant of a maximal order then $D$ is odd or
$D\equiv 8$ or $12\pmod{16}$.
\end{rem}

For $x, y \in \Q$, we write $x \equiv y \pmod{2^m}$ to mean $\ord_2(x-y) \ge m$.

\begin{lem}
\label{trs}
Suppose $\O$ is an order of discriminant $D$ in an imaginary quadratic field $K$, 
$E$ is an elliptic curve over $\C$, 
and $\End(E)=\O$.
Then:
\begin{enumerate}
\item 
there is a 
$\tau \in \H \cap K$ such that $j(\tau) = j(E)$ 
and $\tau = r \td + s$ with 
$r, s \in \Q$, $r \equiv 1 \pmod{2}$, and $s \equiv 0 \pmod{4}$;
\item
with $\tau$ as in {\rm (i)}, then for every 
$\lambda \in \O_{\tau,2}^\times$ we have
$$
\etwo(\lambda) = 
\begin{cases}
\etd(\lambda) & \text{if $r \equiv 1\pmod{4}$,} \\[4pt]
{\etd(\lambda)}(-1)^{(\N_{K/\Q}(\lambda)-1)/2} & \text{if $r \equiv -1\pmod{4}$ 
and $D\equiv 4, 8 \pmod{16}$,} \\[4pt]
\shortoverline{\etd(\lambda)} & \text{if $r \equiv -1\pmod{4}$ and $D \not\equiv 4, 8 \pmod{16}$}.
\end{cases}
$$
\end{enumerate}
\end{lem}

\begin{proof}
By the theory of complex multiplication there is an
invertible ideal $\a \subseteq \O$ such that $j(E) = j(\a)$.  
Changing $\a$ in its ideal class if 
necessary, we may assume that $[\O:\a]$ is odd.  Let $a$ be the smallest 
positive integer in $\a$. Then $\a$ has a $\Z$-basis $\{a,b\td+c\}$ 
with $a,b,c \in \Z$ and $b\td+c\in\H$, and $a, b$ must both be odd.  
Subtracting $ca^2$ from $c$ if necessary, we may assume that $4 \mid c$.  If we let 
$\tau = (b/a)\td + (c/a) \in \H \cap K$ then $L_\tau = a^{-1}\a$, 
so $j(\tau) = j(L_\tau) = j(\a) = j(E)$. This gives (i).
Since $j(\tau) = j(E)$, it follows that $\O_\tau=\O$ ($=\O_{\td}$).

By definition of $q_\tau$, for every $\lambda \in K^\times$ we have 
$$
q_\tau(\lambda) = \tinymatrix{r}{s}{0}{1} q_{\td}(\lambda) \tinymatrix{r}{s}{0}{1}^{-1}.
$$
By Definition \ref{phidef} and the fact that $s \equiv 0 \pmod{4}$, 
if $\lambda \in \O_{\tau,2}^\times$ 
then
\begin{align*}
\eone(\lambda) &= \phi\bigl(\tinymatrix{1}{0}{0}{\N_{K/\Q}(\lambda)^{-1}}
  \tinymatrix{r}{0}{0}{1} q_{\td}(\lambda)\tinymatrix{r}{0}{0}{1}^{-1}\bigr) \\
    &= \phi\bigl(\tinymatrix{r}{0}{0}{1}
    \tinymatrix{1}{0}{0}{\N_{K/\Q}(\lambda)^{-1}}
     q_{\td}(\lambda)\tinymatrix{r}{0}{0}{1}^{-1}\bigr).
\end{align*}
Thus 
$\delta_\tau(\lambda) =\delta_{\tau_D}(\lambda)$
if $r \equiv 1\pmod{4}$, and 
applying Lemma \ref{lemphi}(ii) with
$M =  \tinymatrix{1}{0}{0}{\N_{K/\Q}(\lambda)^{-1}} q_{\td}(\lambda)$ shows 
that $\delta_\tau(\lambda) = \shortoverline{\delta_{\tau_D}(\lambda)}$
if $r \equiv -1\pmod{4}$.
Part (ii)  now follows from 
Definition \ref{e2def} (and the fact that $\etd(\lambda)\in\{\pm 1\}$ when $D\equiv 4, 8 \pmod{16}$).
\end{proof}

\section{$\Q$-curves}
\label{Qc}

Suppose now that $D$ is a (negative) fundamental discriminant, and let 
$d\in\Z^+$ be given by \eqref{ddef} and 
$\td$ by \eqref{zddef}.  Then $d$ is a squarefree positive integer.
With $K := \Q(\sqrt{-d})$, then $\O_{\td} = \O_K$ is the 
maximal order of $K$, and $H := H_{\td}$ is the Hilbert class field of $K$.
Following Gross (\S11 of \cite{grossbk}), an elliptic curve $E$ 
over $H$ is defined to be a {\em $\Q$-curve} if $E$ is isogenous over $H$ to 
$E^\sigma$ for all $\sigma \in \Gal(H/\Q)$.  
By Lemma 11.1.1 of \cite{grossbk}, 
$E$ is a $\Q$-curve if and only if for all but finitely many primes $\P$ 
of $H$ and all $\sigma \in \Gal(H/\Q)$, 
\begin{equation}
\label{qc}
\psi_E(\P^\sigma) = \psi_E(\P)^\sigma
\end{equation}
where $\psi_E$ is the Hecke character of $E$ over $H$.  
In Theorem \ref{qcurvethm} below we use Theorem \ref{main} to exhibit,
whenever $d \equiv 2$ or $3 \pmod{4}$, explicit models and Hecke characters of 
$\Q$-curves, defined over $\Q(j)$, with CM by $\O_K$.
When $d$ is a prime congruent to $3 \pmod{4}$,
Theorem \ref{qcurvethm} was proved by Gross (Theorem 12.2.1 of \cite{grossbk} and 
Proposition 3.5 of \cite{gross}), and when $3 \nmid d \equiv 3 \pmod{4}$ 
it was proved by Stark (Theorem 1 of \cite{Stark}) 
(see Remark \ref{GSrem} below).

\begin{rem}
When all prime divisors of $d > 1$ are congruent to $1\pmod{4}$, 
there are no $\Q$-curves with CM 
by $\O_K$.  See Example 3 on p.~527 of
\cite{shimura1971} and \S11.3 of  \cite{grossbk}.
\end{rem}

We first need a lemma that we will use to prove Theorem \ref{qcurvethm}.

\begin{defn}
If $F$ is a number field, $\q$ is a prime
of $F$, and $a, b \in F^\times$, let 
$\hilb{a}{b}_{\q,F} \in \{\pm1\}$ 
denote the local Hilbert symbol at $\q$, which is defined to be $1$ if and only if 
$b \in \N_{F_\q(\sqrt{a})/F_\q}(F_\q(\sqrt{a})^\times)$. 
Let $\hilbt{a}{b}{F} = \prod_{\q \mid 2}\hilb{a}{b}_{\q,F}$.
\end{defn}

\begin{lem}
\label{homom}
\begin{enumerate}
\item
The function $\etd : \O_{K,2}^\times \to \bmu_4$ 
is a homomorphism.
\item
If $d \equiv 3 \pmod{4}$ and $\lambda \in \O_K$ is prime to $2$, 
then 
$$
\etd(\lambda) = \hilbt{\sqrt{-d}}{\lambda}{K}.\phantom{(-1)^{(a-1)/2}}\phantom{(-1)^{(q-1)(q+d+3)/16}}
$$
\item
If $d \equiv 6 \pmod{8}$, $\lambda \in \O_K$ is prime to $2$, 
and $q = \N_{K/\Q}(\lambda)$, then
$$
\etd(\lambda) = (-1)^{(q-1)(q+d+11)/16}\hilbt{\sqrt{-d}}{\lambda}{K}.\phantom{(-1)^{(a-)/2}}
$$
\item
If $d \equiv 2 \pmod{8}$, $u, v \in \Z$,
$\lambda = u + v \sqrt{-d}$ is prime to $2$, 
and $q = \N_{K/\Q}(\lambda)$, then
$$
\etd(\lambda) = (-1)^{(u-1)/2}(-1)^{(q-1)(q+d+3)/16}\hilbt{\sqrt{-d}}{\lambda}{K}.
$$
\end{enumerate}
\end{lem}

\begin{proof}
Part (i) can be checked directly using Proposition \ref{eprop}.
It is easy to check that both sides of the displayed equations
depend only on $\lambda\pmod{8\O_K}$, so 
(ii), (iii), and (iv) can also be checked by direct computations.
\end{proof}

Let $j=j(\td)$, $\gamma_2 = \gamma_2(\td)$, and $\gamma_3 = \gamma_3(\td)$.

\begin{thm}
\label{qcurvethm}
Suppose $d \equiv 2$ or $3 \pmod{4}$.  Let $E$ be the curve
$$
E = 
\begin{cases}
\Et{\sqrt{-d}\gammatwod^4}{\td} : y^2 = x^3 + \frac{d {\jd}^3}{48} x 
    - \frac{d\sqrt{-d}{\gammathreed}{\jd}^4}{864} & \text{if $d \equiv 3 \pmod{4}$,} \\
\Et{-\sqrt{d}\gammatwod^4}{\td} : y^2 = x^3 - \frac{d {\jd}^3}{48} x 
    - \frac{d\sqrt{d}{\gammathreed}{\jd}^4}{864} & \text{if $d \equiv 2 \pmod{4}$.}
    \end{cases}
$$
Then:
\begin{enumerate}
\item
$E$ is defined over $\Q({\jd})$.
\item
\label{jDeltad}
$j(E) = {\jd}$ and 
$\Delta(E) = (-1)^dd^3{\jd}^8$.
\item
$E$ is a $\Q$-curve.
\item
\label{psiqpart}
Suppose $\P$ is a prime of $H$, not dividing $2$, where $E$ has good reduction.
Suppose
$\lambda = u+v \sqrt{-d} \in \O_K$ is a generator of $\N_{H/K}(\P)$, with $u, v \in \frac{1}{2}\Z$,
and let $q = \N_{H/\Q}(\P) = u^2 + d v^2$.  If $d \neq 3$ then the Hecke character $\psi$ 
of $E$ over $H$ is given by
$$
\psi(\P) = 
\begin{cases}
\leg{4u}{d}\lambda & \text{if $d \equiv 3 \pmod{4}$} \\
(-1)^{(q-1)(q+d+11)/16}\leg{u}{d/2}\lambda & \text{if $d \equiv 6 \pmod{8}$} \\
(-1)^{(u-1)/2}(-1)^{(q-1)(q+d+3)/16}\leg{u}{d/2}\lambda & \text{if $d \equiv 2 \pmod{8}$}
\end{cases}
$$
where $\leg{~}{~}$ is the Jacobi symbol.
\end{enumerate}
\end{thm}

\begin{proof}
Note that $E$ is the curve of Corollary \ref{exrem} with $\beta = \sqrt{-d}$.  
By Lemma \ref{notjlem}(ii,iii) we have (i). 
By \eqref{jEalph} and \eqref{g23a} we have (ii).

Suppose $\psi$, $\P$, $\lambda$, $q$ and $u$ 
are as in \ref{psiqpart}.
By Corollary \ref{exrem}(ii) (with $\beta = \sqrt{-d}$), 
\begin{equation}
\label{gr-1}
\psi(\P) = \legtF{\sqrt{-d}}{\P}{H}\etd(\lambda)\lambda.
\end{equation}
We will evaluate $\legt{\sqrt{-d}}{\P}$ using quadratic reciprocity over $K$.

Let $\p$ be the prime of $K$ below $\P$ and let $f = [\O_H/\P:\O_K/\p]$, so 
$\lambda\O_K = \N_{H/K}(\P) = \p^f$.  
By Proposition II.7.4.3(v,viii) of \cite{Gras} and the product formula, 
\begin{equation}
\label{gr0}
\legtF{\sqrt{-d}}{\P}{H} = 
\bigl(\textstyle\frac{\sqrt{-d}}{\p}\bigr)^f_{\hskip-1pt 2,K}
= \hilb{\sqrt{-d}}{\lambda}_{\p,K}
= \prod_{\q \ne \p}\hilb{\sqrt{-d}}{\lambda}_{\q,K}
\end{equation}
where $\q$ runs over primes of $K$.
If $\q \nmid 2d$ then $\q$ is unramified in $K((\sqrt{-d})^{1/2})/K$.
Since $\ord_\q(\lambda)=0$ for all $\q \neq \p$, it follows from
Proposition II.7.1.1(vi) of \cite{Gras} that if $\q \nmid 2 \p d$ then
$\hilb{\sqrt{-d}}{\lambda}_{\q,K} = 1$, so
\begin{equation}
\label{gr1}
\prod_{\q \ne \p}\hilb{\sqrt{-d}}{\lambda}_{\q,K}
    = \prod_{\q \mid d, \q \nmid 2}\hilb{\sqrt{-d}}{\lambda}_{\q,K}
      \prod_{\q \mid 2}\hilb{\sqrt{-d}}{\lambda}_{\q,K}.
\end{equation}
Suppose $\q \mid d$ and $\q \nmid 2$.  Then 
$\lambda \equiv u \pmod{\q\O_{K_\q}}$ and 
$\hilb{\sqrt{-d}}{\lambda}_{\q,K}  = \hilb{\sqrt{-d}}{u}_{\q,K}$.
Further, $\q$ ramifies in $K/\Q$, so if $\ell = \N_{K/\Q}(\q)$, then
$$
\hilb{\sqrt{-d}}{u}_{\q,K}  = \hilb{d}{u}_{\ell,\Q} = \leg{4u}{\ell},
$$
the first equality by Proposition II.7.1.1(ii,iv) of \cite{Gras}, and the second 
by Theorem 1 in \S{III.1.2} of \cite{serre} (and the fact that $u$ is a half-integer).
Thus if $d'$ is the largest odd divisor of $d$ 
and $\ell$ runs over primes of $\Q$, then 
\eqref{gr0} and \eqref{gr1} yield
$$
\legtF{\sqrt{-d}}{\P}{H} 
    = \prod_{\ell \mid d'}\leg{4u}{\ell} 
        \prod_{\q \mid 2}\hilb{\sqrt{-d}}{\lambda}_{\q,K}
    = \leg{4u}{d'}\hilbt{\sqrt{-d}}{\lambda}{K}.
$$
Combining this with \eqref{gr-1} 
gives
$$
\psi(\P) = \leg{4u}{d'}\hilbt{\sqrt{-d}}{\lambda}{K} \etd(\lambda) \lambda.
$$
Now \ref{psiqpart} follows from Lemma \ref{homom}.

To prove that $E$ is a $\Q$-curve, we need to check that \eqref{qc} holds
for all primes $\P$ of $H$ as above and all $\sigma \in \Gal(H/K)$.  
This is clear from the formulas of \ref{psiqpart}. 
\end{proof}

By Proposition \ref{Eorder}\ref{Eorderpts}, Theorem \ref{qcurvethm}\ref{psiqpart} 
gives formulas for $|E(\O_K/\P)|$.

\begin{rem}
\label{GSrem}
Suppose that $d \equiv 2$ or $3 \pmod{4}$, and suppose that $3 \nmid d$.  
Let $A$ be the elliptic curve
$$
A = 
\begin{cases}
\Et{\sqrt{-d}}{\td} : y^2 = x^3 + \frac{d{\gammatwod}}{48}x 
    - \frac{d \sqrt{-d} {\gammathreed}}{864} & \text{if $d \equiv 3 \pmod{4}$,} \\
\Et{-\sqrt{d}}{\td} : y^2 = x^3 - \frac{d{\gammatwod}}{48}x 
    - \frac{d\sqrt{d} {\gammathreed}}{864} & \text{if $d \equiv 2 \pmod{4}$.}
\end{cases}
$$
By \S6 of \cite{Birch} or Theorem 2 of \cite{Schertz}, ${\gammatwod} \in \Q({\jd})$ 
(this is where $3 \nmid d$ is used), 
so $A$ is defined over $\Q({\jd})$ and is  isomorphic over $\Q({\jd})$ to 
the $E$ of Theorem \ref{qcurvethm}.  
By \eqref{jEalph} and Lemma \ref {notjlem}(i), $j(A) = {\jd}$ and 
$\Delta(A) = -d^3$, and $A$ is a $\Q$-curve by Theorem \ref{qcurvethm}(iii).  
When $d$ is a prime $p$, $A$ is the model given 
by Gross in \cite{grossbk,gross} for the $\Q$-curve that he denoted $A(p)$.  
When $3 \nmid d$ and $d \equiv 7 \pmod{8}$ (resp., $d \equiv 3 \pmod{8}$), 
$A$ is the curve $E_1$ (resp., $E_{-1}$) considered by Stark in 
Theorem 1 of \cite{Stark}.
\end{rem}

\section{Elliptic curves over $\Fp$ with $p \equiv 1 \pmod{4}$}
\label{yamihere}

Theorem \ref{1mod4thm} below, which
uses Theorem \ref{main}, gives a simple formula for the number of points on
an ordinary elliptic curve $E$ over $\Fp$ when $p \equiv 1 \pmod{4}$ and
$\End_{\bar\F_p}(E)=\O_{\Q({\sqrt{-d})}}$ with 
$d \equiv 2$ or $3 \pmod{4}$.

We will use the following lemma, which is a variant of Deuring's Lifting Theorem.

\begin{lem}
\label{ugly}
Suppose $p$ is prime, $E$ is an ordinary elliptic curve over $\Fp$, and
$\O := \End_{\Fp}(E)$ is an order in an imaginary quadratic
field $K$.  Let $H=K(j(\O))$.
Then there are an elliptic curve $\E$ over $H$ and a prime $\P$ of $H$ such that
$\O_H/\P \cong \F_p$, 
$\End_{H}(\E)= \O$, $j(\E) = j(\O)$, and 
the reduction of $\E$ modulo $\P$ is isomorphic to $E$ over $\F_p$.
\end{lem}

\begin{proof}
Since the proof is easy when $j=0$ or $1728$, we can reduce to the case $\O^\times =\{\pm 1\}$.
Since $E$ is ordinary, $E$ has a canonical lifting $\Ec$ to $\Qp$ 
(see Theorem 3.3 on p.~172 of \cite{lifting}), 
i.e., $\Ec$ is an elliptic curve over $\Qp$ that reduces to $E$, 
and $\End_{\Q_p}(\Ec) = \End_{\Fp}(E) = \O$.
The action of $\End_{\Q_p}(\Ec)$ on the space $\Omega$ of holomorphic differentials induces
an embedding $K\cong \End_{\Q_p}(\Ec)\otimes\Q \hookrightarrow \End(\Omega) \cong\Qp$.  
By the theory of complex multiplication (see Theorem 5.7(iii) of \cite{Shimurabook}), 
we can fix an embedding $\Qp \hookrightarrow \C$ under which $j(\Ec) = j(\O)$.
Since $K \subset \Qp$ and $j(\O) = j(\Ec) \in \Qp$,
we have $H = K(j(\O)) \subset \Qp$.  
Let $\P = \O_H \cap p\Z_p$.  Then $\P$ is a prime of $H$ with residue field 
$\O_H/\P \cong \Fp$. 
Since $\Ec$ is a lift of $E$, $j(\O) = j(\Ec)$ reduces to $j(E)$ modulo $\P$.  

Let $\A$ be an elliptic curve over $H$ with $j(\A) = j(\O)$.
Then $\Ec$ is a quadratic twist of $\A$ by some $\delta\in\Q_p^\times$.
Choose ${\delta}' \in\Q^\times$ so that $u:={\delta}'/\delta$ is in $\Z_p^\times$ and
let $\E$ be the quadratic twist of $\A$ by ${\delta}'$.
Then $\Delta(\E)=u^6\Delta(\Ec)$, which is in $\Z_p^\times$
since $\Ec$ has good reduction at $p$.
Thus $\E$ is an elliptic curve over $H$ 
with good reduction at $\P$ and with $j(\E) = j(\O)$.
In particular, $\End_H(\E) = \O$. 
Since the reduction $\tilde\E$ of $\E$ modulo $\P$ 
has $j$-invariant $j(E)$, and $\Aut(E) = \O^\times = \{\pm1\}$, 
it follows that $\tilde\E$ is a quadratic twist 
of $E$.  Thus replacing $\E$ by a quadratic twist ensures that 
$\tilde\E$ is isomorphic to $E$ over $\O_H/\P = \Fp$.  
\end{proof}

If $a \in \Fp^\times$ is a square, let $\legf{a}{p}$ be the quartic residue symbol 
defined by 
$$
\legf{a}{p} \in \{\pm1\}, \quad \legf{a}{p} \equiv a^{(p-1)/4} \pmod{p}.
$$

\begin{thm}
\label{1mod4thm}
Suppose $p$ is prime, $E$ is an ordinary elliptic curve over $\Fp$, and 
$\O := \End_{\Fp}(E)$ is an order in an imaginary quadratic field $K$.  
Suppose further that $p \equiv 1 \pmod{4}$, and the discriminant $D$ of $\O$ 
is either odd and not $-3$, or is congruent to $4$ or $8 \pmod{16}$.  
Then:
\begin{enumerate}
\item
the discriminant $\Delta(E)$ of $E$ is a square in $\Fp^\times$,
\item
there are $u, v \in \frac{1}{2}\Z$ 
such that $u^2 + |D| v^2 = p$ and 
$\lambda := u + v \sqrt{D} \in \O$ satisfies
$$
\renewcommand{\arraystretch}{1.35}
\arraycolsep=1pt
\begin{array}{rlll}
&\lambda^3 &\equiv 1 \pmod{4\O} &\text{if $D$ is odd}, \\
(-1)^{(p-1)/4}&\lambda &\equiv 1 \text{~or $1 + \sqrt{D} \pmod{4\O}$} &\text{if $D \equiv 4 \pmod{16}$}, \\
&\lambda &\equiv 1 \text{~or $-1 + \sqrt{D} \pmod{4\O}$~} &\text{if $D \equiv 8 \pmod{16}$},
\end{array}
$$
\item
if $u$ is as in {\rm (ii)}, then
$
|E(\Fp)| = p+1 - 2\legf{\Delta(E)}{p}u.
$
\end{enumerate}
\end{thm}

\begin{proof} 
Let $j = j(\O)$ and $H = K(j)$.  
Using Lemma \ref{ugly}, fix an elliptic curve $\E : y^2 = x^3 + ax + b$ 
over $H$ and a prime $\P$ of $H$ such that $\O_H/\P \cong \F_p$, $j(\E) = j$, and 
the reduction of $\E$ modulo $\P$ is isomorphic over $\F_p$ to $E$.  
Let $\p = \P \cap K$.  Since $\O_H/\P \cong \F_p$,  we have 
$\N_{H/K}(\P) = \p$, so $\p$ is principal with a generator 
$\lambda = u+v\sqrt{D} \in \O$. 
In particular $u^2 + |D|v^2 = \N_{K/\Q}(\lambda) = p$.
Since $p \equiv 1 \pmod{4}$, we have $\P \nmid 2$.

Suppose first that $D$ is odd.  Then $(\O/2\O)^\times \cong (\O_K/2\O_K)^\times$ has 
order $1$ or $3$, so $\lambda^3 \equiv 1 \pmod{2\O}$.  Further, 
$\N_{K/\Q}(\lambda^3) = p^3 \equiv 1 \pmod{4}$.
A straightforward computation shows that
the only elements  in $(\O/4\O)^\times$ that are $1$ mod $2\O$
and have norm $1$ are $\pm1$, so
$\lambda^3 \equiv \pm1 \pmod{4\O}$.  Replace $\lambda$ by $-\lambda$, 
if necessary, to ensure that $\lambda^3 \equiv 1 \pmod{4\O$}.

Now suppose $D \equiv 4$ or $8 \pmod{16}$.  Since 
$\N_{K/\Q}(\lambda) = p \equiv 1 \pmod{4}$, a straightforward computation 
in $(\O/4\O)^\times$ shows that 
$\lambda \equiv \pm1$ or $\pm1+\sqrt{D} \pmod{4\O}$.  
Replace $\lambda$ by $-\lambda$, 
if necessary, to ensure that $(-1)^{(p-1)/4}\lambda \equiv 1$ or $1+\sqrt{D} \pmod{4\O$} 
when $D \equiv 4 \pmod{16}$, and $\lambda \equiv 1$ or $-1+\sqrt{D} \pmod{4\O$} 
when $D \equiv 8 \pmod{16}$.  
Note that if $D \equiv 8 \pmod{16}$ then $p \equiv 1 \pmod{8}$.

Thus we have (ii).  Note that 
if $u', v' \in \frac{1}{2}\Z$ is another pair satisfying (ii), then 
$u'+v'\sqrt{D} \in \O$ is a generator of a prime of $K$ above $p$, 
so $u' = \pm u$ and $v' = \pm v$. By the 
congruences on $\lambda$ in (ii), we have $u' = u$, i.e., the $u$ 
satisfying (ii) is unique.

Let $\td$ be as defined by \eqref{zddef}.  
We will apply Theorem \ref{main} to $\E$ with $\tau = \td$.  
Let $\nu = 1$ if $D$ is odd, and 
$\nu = i$ if $D$ is even.  
By Proposition \ref{eprop}, $\etd(\lambda) = \nu^{(p-1)/2}$ 
(we use here that $p \equiv 1 \pmod{8}$ if $D \equiv 8 \pmod{16}$).
By Theorem \ref{main}, since $\Tr_{K/\Q}(\lambda)=2u$,
\begin{equation}
\label{le}
|E(\Fp)| = |\E(\O_H/\P)| 
   = p+1-2\nu^{(p-1)/2}\legt{6b\nu{\gammathreed}}{\P}u.
\end{equation}
Note that
$(2^53^3b)^2/\Delta(\E) = j(\E)-1728= \pm(\nu{\gammathreed})^2$.
It follows from Lemma \ref{notjlem}(ii,iii) and $p \equiv 1 \pmod{4}$ that 
modulo $\P$, $\Delta(\E)$ is a square and
\begin{multline*}
\nu^{(p-1)/2}\legt{6b\nu{\gammathreed}}{\P} 
    \equiv \nu^{(p-1)/2}(6b\nu{\gammathreed})^{(p-1)/2} 
    \equiv (6^2 b^2 (j(\E)-1728))^{(p-1)/4} \\
    = (2^{12}3^8 b^4/\Delta(\E))^{(p-1)/4} 
    \equiv \legf{\Delta(E)^{-1}}{p} = \legf{\Delta(E)}{p}.
\end{multline*}
Since the outer terms are $\pm1$, they must be equal.  
Now combine this with \eqref{le}.  
\end{proof}

\begin{rem}
With notation as in Theorem \ref{1mod4thm}, if $E$ is supersingular 
rather than ordinary, and if further $p \ge 5$, then $|E(\Fp)| = p+1$.
\end{rem}

\end{document}